\input amstex
\documentstyle{amsppt}
\NoBlackBoxes
\magnification1200
\pagewidth{6.5 true in}
\pageheight{9 true in}

\topmatter
\title
Real zeros of  quadratic Dirichlet $L$-functions
\endtitle
\author
J. B. Conrey and K. Soundararajan
\endauthor
\address
American Institute of Mathematics, 360 Portage Avenue, Palo Alto CA , USA
\endaddress
\email
conrey\@ best.com
\endemail
\address
Department of Mathematics, University of Michigan, Ann Arbor, MI 48109, USA
\endaddress
\email
ksound\@ math.ias.edu
\endemail
\thanks Research of both authors is supported
by the American Institute of Mathematics (AIM), and in part
through grants from the NSF including FRG grant DMS-00-74028.
\endthanks
\endtopmatter

\def\sumflat{\sideset \and^{\flat} \to \sum}
\def\L{\fracwithdelims()}

\def\phi{\varphi}
\def\lam{\lambda}

\def\Sum{\Cal S}

\head 1. Introduction \endhead

\noindent  A small part of the Generalized Riemann Hypothesis
asserts that $L$-functions do not have zeros on the
line segment $(\frac 12 ,1]$.  The question of vanishing at $s=\frac 12$
often has deep arithmetical significance, and
has been investigated extensively.  A persuasive view is
that $L$-functions vanish at $\frac 12$ either for trivial
reasons (the sign of the functional equation being negative), or
for deep arithmetical reasons (such as the $L$-function of
an elliptic curve of positive rank) and that the latter
case happens very rarely.  N. Katz and P. Sarnak [7] have
formulated precise conjectures on the low lying zeros
in families of $L$-functions which support this view.

In the case of Dirichlet $L$-functions it is expected that
$L(\tfrac 12,\chi)$ is never zero, and so $L(\sigma,\chi)\neq 0$ for
all $\tfrac 12 \le \sigma \le 1$.  This conjecture appears to
have been first enunciated by S.D. Chowla [2] in the special
case of quadratic characters $\chi$.  Progress towards these
non-vanishing questions has been in two directions:
zero-density type results which establish that very few $L$-functions
have a zero in $(\frac 12 +\epsilon,1]$ (see for example
A. Selberg [10], M. Jutila [6] and D.R. Heath-Brown [4]), and
a growing body of work on non-vanishing at $\frac 12$
(see for example
R. Balasubramanian and V.K. Murty [1], H. Iwaniec and Sarnak [5],
and K. Soundararajan [11]).  Further much numerical
evidence for the GRH has been accumulated, and these
calculations support Chowla's conjecture (see [8] and [9]).
However the state of knowledge could not exclude the possibility
that every Dirichlet $L$-function of sufficiently large conductor
has a non-trivial real zero.  In this paper we eliminate
this possibility, and prove that a positive proportion of
quadratic Dirichlet $L$-functions do not vanish on $[\frac 12,1]$.

For an integer $d\equiv 0$, or $1 \pmod 4$ we put $\chi_d
(n)=\L{d}{n}$, so
that $\chi_d$ is a real character with conductor at most $|d|$.
If $d$  is an odd, positive, square-free integer
then $\chi_{-8d}$ is a real, primitive character with conductor
$8d$, and with $\chi_{-8d}(-1)=-1$.

\proclaim{Theorem 1}  For at least $20\%$ of the odd
square-free integers $d\ge 0$ we have $L(\sigma,\chi_{-8d}) >0$
for $0\le \sigma \le 1$.  More precisely, for all large $x$ the
number of odd positive square-free integers $d \le x$
such that $L(\sigma,\chi_{-8d}) >0$ for all $0\le \sigma \le 1$
exceeds $\frac{1}{5} (\frac {4x}{\pi^2})$.
\endproclaim

While in this paper we have restricted our attention
to fundamental discriminants of the form $-8d$, our methods would
apply to fundamental discriminants in any arithmetic progression.
Also our proof yields that there are many $L$-functions
having no non-trivial zeros in a thin rectangle containing the
real axis.  Precisely, there is a constant $c>0$ such that
for at least $20\%$ of the fundamental discriminants $-8d$ with $0< d\le x$,
the rectangle $\{\sigma+it: \ \ \sigma \in [0,1], \ \ |t|\le c/\log x\}$
is free of zeros of $L(s,\chi_{-8d})$.
As another consequence of our work we find that the number of
fundamental discriminants $-8d$ with $0 <d\le x$ such
that $L(s,\chi_{-8d})$ has a zero in the
interval $[\sigma,1]$ is $\ll x^{1-(1-\epsilon)(\sigma-\frac 12)}$
for any fixed $\epsilon >0$.

\head 2.  Outline of the Proof \endhead

\noindent We begin with the following version of
the argument principle, due to Selberg [10], whose
proof we reproduce for completeness.

\proclaim{Lemma 2.1}  Let $f$ be a holomorphic function, which is
non-zero in some half-plane Re$(z) \ge W$.  Let
${\Cal B}$ be the rectangular box with vertices $W_0\pm iH$,
$W_1 \pm iH$ where $H>0$ and $W_0 < W < W_1$.  Then
$$
\align
4H
\sum\Sb \beta +i\gamma \in {\Cal B}\\ f(\beta+i\gamma) =0\endSb
&\cos \L{\pi\gamma}{2H} \sinh \L{\pi(\beta-W_0)}{2H}
=\int_{-H}^{H} \cos \L{\pi t}{2H} \log |f(W_0+it)| dt
\\
&+\int_{W_0}^{W_1}\sinh \L{\pi (\alpha-W_0)}{2H}
\log |f(\alpha+iH)f(\alpha-iH)| d\alpha
\\
&-\text{Re} \int_{-H}^{H} \cos \Big(\pi \frac{W_1-W_0 +it}{2iH}\Big)
\log f(W_1+it) dt.\\
\endalign
$$
\endproclaim
\demo{Proof}  From the box ${\Cal B}$ we exclude the line segments
$x+i\gamma$ with $W_0\le x\le \beta$ for every zero $\beta+i\gamma$
of $f$ lying in ${\Cal B}$.  Denoting by $\Gamma$ the
boundary of the resulting domain we see that
$$
\int_{\Gamma} \cos \Big( \pi \frac{s-W_0}{2iH} \Big) \log f(s) ds
= 0.
$$
Since the value of $\log f(s)$ differs by $2\pi i$ on the upper and lower
sides of the ``cuts'' from $\beta+i\gamma$ to $W_0+i\gamma$, we
conclude from the above that
$$
\align
2\pi i \sum\Sb \beta +i\gamma \in {\Cal B}\\ f(\beta+i\gamma) =0\endSb
&\int_{W_0+i\gamma}^{\beta+i\gamma} \cos \Big( \pi \frac{s-W_0}{2iH}
\Big) ds \\
= \biggl(& \int_{W_0-iH}^{W_0+iH} + \int_{W_0+iH}^{W_1+iH} -
\int_{W_1-iH}^{W_1+iH} - \int_{W_0-iH}^{W_1-iH}\biggr)
 \cos \Big( \pi \frac{s-W_0}{2iH}\Big) \log f(s) ds.
\\
\endalign
$$

The imaginary part of the LHS above equals the LHS of the
equality of the Lemma.
%of (6.1) is
%$$
%\align
%&\text{Im} \biggl( 4H
%\sum\Sb \beta +i\gamma \in {\Cal B}\\ f(\beta+i\gamma) =0\endSb
%\biggl( \cos \Big(\pi \frac{\beta-W_0}{2iH} + \pi \frac{H+\gamma}{2H}\Big)
%-\cos \Big(\pi\frac{H+\gamma}{2H}\Big)\biggr)\biggr) \\
%&= 4H
%\sum\Sb \beta +i\gamma \in {\Cal B}\\ f(\beta+i\gamma) =0\endSb
%\cos \L{\pi\gamma}{2H} \sinh \L{\pi(\beta-W_0)}{2H}.\\
%\endalign
%$$
The imaginary part of the first integral on the RHS above
equals the first term on the RHS of the Lemma.
% is
%$$
%\text{Re} \int_{-H}^{H} cos \L{\pi t}{2H} \log f(W_0+it) dt =
%\int_{-H}^{H} \cos \L{\pi t}{2H} \log |f(W_0+it)| dt.
%$$
The second and fourth integrals on the RHS above have combined imaginary part
equal to the second term on the RHS of the Lemma.
%$$
%\text{Im} \int_{W_0}^{W_1} \sin\Big( \pi \frac{\alpha-W_0}{2iH} + \pi\Big)
%\log f(\alpha+iH) d\alpha
%= \text{Im}  \biggl( \int_{W_0}^{W_1} i\sinh \L{\pi (\alpha-W_0)}{2H} \log
%f(\alpha+iH) d\alpha \biggr)
%\\
%= \int_{W_0}^{W_1}\sinh \L{\pi (\alpha-W_0)}{2H} \log |f(\alpha+iH)| d\alpha.
%$$
Lastly the imaginary part of the fourth term on the RHS above equals
the third term on the RHS of the Lemma.
%$$
%\text{Im} \biggl(- \int_{W_0}^{W_1} \sin\Big( \pi \frac{\alpha-W_0}{2iH}\Big)
%\log f(\alpha-iH) d\alpha \biggr) =
%\int_{W_0}^{W_1}\sinh \L{\pi (\alpha-W_0)}{2H} \log |f(\alpha-iH)| d\alpha.
%$$
%The third term in the RHS of (6.1) has imaginary part
%$$
%-\text{Re} \int_{-H}^{H} \sin\Big(\pi \frac{W_1-W_0 +i(H+t)}{2iH}\Big)
%\log f(W_1+it) dt.
%$$
Thus Lemma 2.1 is proved.

\enddemo

Let $X$ be large, and let
$d$ be any odd square-free number in $[X,2X]$.  We shall apply
Lemma 2.1 to a mollified version of $L(s,\chi_{-8d})$.  Precisely,  for
a parameter $X^{\epsilon} \le M \le X$ to be fixed later
\footnote {Here and throughout, $\epsilon$ denotes a small positive real
number. The reader should be warned that it might be a different $\epsilon$
from line to line.},
let
$$
M(s,d) = \sum_{n\le M} \frac{\lam(n)}{n^s} \chi_{-8d}(n),
$$
where the $\lambda(n)$ are real numbers $\ll n^{\epsilon}$ to be specified
later.  We apply Lemma 2.1 with
$f(s,d) := L(s,\chi_{8d}) M(s,d)$ and $W_0=\frac 12 -\frac{R}{\log X}$,
$H=\frac{S}{\log X}$, and $W_1=\sigma_0$ where $R$ and $S$ are
fixed positive parameters in the interval $(\epsilon, 1/\epsilon)$
to be chosen later, and $\sigma_0 >1$ will
be chosen later such that $f(s,d)$ has no zeros in Re $ s> \sigma_0$.
Since the LHS of Lemma 2.1 consists of positive terms
we glean that
$$
4S \sum\Sb \beta\ge \frac 12 -\frac{R}{\log X} \\
L(\beta,\chi_{-8d})=0\endSb \sinh \Big( \frac{\pi(R+\log X(\beta-1/2))}{2S}
\Big) \le I_1(d) + I_2(d) +I_3(d), \tag{2.1}
$$
where (after obvious changes of variables)
$$
I_1(d) =
\int_{-S}^{S} \cos \L{\pi t}{2S} \log \Big|f\Big(\frac 12 -\frac{R}{\log X}
+ i\frac{t}{\log X},d\Big)\Big|dt, \tag{2.2a}
$$
$$
I_2(d) = \int_{-R}^{(\sigma_0-\frac 12) \log X} \sinh
\Big(\frac{\pi(x +R)}{2S}\Big)
\log \Big|f\Big(\frac 12+\frac{x}{\log X}+i\frac{S}{\log X},d\Big)
\Big|^2 dx, \tag{2.2b}
$$
and
$$
I_3(d) = -\text{Re }
\int_{-S}^{S} \cos\Big(\pi\frac{(\sigma_0-1/2)\log X -R +it}{2iS}\Big)
\log f\Big(\sigma_0 +i\frac t{\log X},d\Big) dt.
\tag{2.2c}
$$

Suppose that $L(\beta,\chi_{-8d}) =0$ for some $\beta \in [\frac 12,1]$.
We claim that the LHS of (2.1) exceeds $8S\sinh (\frac{\pi R}{2S})$.
To see this, suppose first that $L(s,\chi_{-8d})$ has a zero
$\beta > \frac 12 +\frac R{\log X}$.  Then the contribution of this zero
alone would be $\ge 4S \sinh \L {\pi R}{S} \ge 8S\sinh \L{\pi R}{2S}$ since
$\sinh(2x) \ge 2\sinh x$ for $x\ge 0$.  On the other hand,
if $L(s,\chi_{-8d})$ has a zero at $\frac 12 +\frac{\xi}{\log X}$
for some $0\le \xi\le R$ then by the functional equation it
also has a zero at $\frac 12 -\frac{\xi}{\log X}$.
In case $\xi =0$ note that there is at least a
double zero at $\frac 12$.  Both these zeros are included in the
LHS of (2.1), and together they contribute
$4S \Big( \sinh \L{\pi(R-\xi)}{2S} +\sinh\L{\pi(R+\xi)}{2S}
\Big) \ge 8S\sinh \L{\pi R}{2S}$,
since the minimum value of $\sinh(x-y) +\sinh(x+y)$ for
$0\le y\le x$ is attained at $y=0$.  We document this below:
$$
I_1(d)+I_2(d)+I_3(d) \ge 8S \sinh\L{\pi R}{2S} \qquad \text{if }
L(s,\chi_{-8d}) \text{  has a non-trivial real zero}. \tag{2.3}
$$

The plan now is to obtain upper bounds for $I_1(d)+I_2(d)+I_3(d)$
on average over $d$, and thereby conclude that the inequality (2.3)
cannot hold too often.  To elaborate on this, we first fix some notation.
Let $\{a_n\}_{n=1}^{\infty} $ be any sequence of complex numbers,
and let $F$ denote a smooth function supported in the interval $[1,2]$.
Throughout this paper we adopt the notation
$$
{\Cal S}(a_d;F) = {\Cal S}(a_d;F,X) = \frac{1}{X} \sum_{d \ \text{odd} }
\mu^2(d) a_d F\L{d}{X}.
$$
Let $\Phi$ be a smooth non-negative function supported in $[1,2]$.
For a complex number $w$ we define
$$
{\check \Phi}(w) = \int_{0}^{\infty} \Phi(y)y^w dy. \tag{2.4a}
$$
For integers $\nu \ge 0$ we define
$$
\Phi_{(\nu)} = \max_{0\le j\le \nu} \int_1^2 |\Phi^{(j)}(t)|dt. \tag{2.4b}
$$
Integrating by parts $\nu$ times we get that
$$
{\check \Phi}(w) =\frac{1}{(w+1)\cdots (w+\nu)} \int_0^{\infty}
\Phi^{(\nu)}(y)y^{w+\nu} dy,
$$
so that for Re $w > -1$ we have
$$
|{\check \Phi}(w)| \ll \frac{2^{\text{Re }w}}{|w+1|^\nu} \Phi_{(\nu)}.
\tag{2.4c}
$$

Let ${\Cal N}(X,\Phi)$ count, with weight $\Phi(d/X)$, the
odd,  square-free integers $d\in [X,2X]$ such that
$L(s,\chi_{-8d})$ has a non-trivial real zero.
In view of (2.1) and (2.3) we see that
$$
{\Cal N}(X,\Phi) \le X\Big( 8S \sinh \Big(\frac{\pi R}{2S}\Big)\Big)^{-1}
{\Cal S} (I_1(d)+I_2(d)+I_3(d);\Phi).
$$
For a complex number $\delta_1$ we define
$$
{\Cal W}(\delta_1,\Phi) =
\frac{{\Cal S}(|L(\tfrac 12+\delta_1,\chi_{-8d})M(\tfrac 12
+\delta_1,\chi_{-8d})|^2;\Phi)}{{\Cal S}(1;\Phi)}. \tag{2.5}
$$
Since the arithmetic mean exceeds the geometric mean we have that
$$
{\Cal S}(\log |f(\tfrac 12+\delta_1,d)|^2;\Phi) \le {\Cal S}(1;\Phi)
\log {\Cal W}(\delta_1;\Phi).
$$
Using this in (2.1), and recalling the definitions (2.2a,b), we
conclude that
$$
{\Cal N}(X,\Phi) \le
\frac{X{\Cal S}(1;\Phi)}{8S \sinh (\frac{\pi R}{2S})}
\Big( J_1(X;\Phi) + J_2(X;\Phi)\Big) + \frac{X}{8S \sinh(\frac{\pi R}{2S}
)} {\Cal S}(I_3(d);\Phi), \tag{2.6}
$$
where
$$
J_1(X;\Phi) = \int_0^S \cos \L{\pi t}{2S}
\log {\Cal W}
\Big(-\frac{R}{\log X} + i\frac{t}{\log X};\Phi\Big) dt, \tag{2.7a}
$$
$$
J_2(X;\Phi) = \int_{-R}^{(\frac 12+\epsilon)\log X}
\sinh \Big( \frac{\pi (x+R)}{2S}\Big)
\log {\Cal W}\Big(\frac{x}{\log X} + i\frac{S}{\log X};\Phi\Big) dx. \tag{2.7b}
$$

At this juncture we specify more carefully the
choice of our mollifier coefficients.  To counter
the rapid growth of the $\sinh (\pi(x+R)/(2S))$ term in (2.7b), we
would like ${\Cal W}((x+iS)/\log X;\Phi)$ to tend rapidly to $1$.
One way to ensure this is to choose $\lam(n)=0$ if $n$ is even, or if $n>M$,
and for odd integers $n \le M$ define
$$
\lam(n) :=
\mu(n) Q\Big(\frac{\log (M/n)}{\log M}\Big)
:=
\cases
\mu(n) & \text {if  } n\le M^{1-b}\\
\mu(n) P(\frac{\log (M/n)}{\log M}) &\text{if } M^{1-b} \le n\le M.\\
\endcases
\tag{2.8}
$$
Here $b$ is a parameter in $[\epsilon, 1-\epsilon]$, and
$P(x)$ is a polynomial such that $P(0)=P^{\prime}(0)=0$, and
$P(b)=1$, $P^{\prime}(b)=0$.

\proclaim{Proposition 2.1}  Suppose $\Phi$ is a non-negative smooth function
supported on $[1,2]$ such that $\Phi(t) \ll 1$, and
with $\int_1^2 \Phi(t) dt \gg 1$.  If $M\le \sqrt{X}$ and
$\delta_1$ is a complex number with $\text{Re }\delta_1 >\epsilon$ then
$$
{\Cal W}(\delta_1,\Phi) = 1 + O(\Phi_{(2)} X^{\epsilon}
(M^{-2\text{Re }\delta_1 (1-b)}+ M^{(\frac 12-\delta_1)(1-b)}X^{-\frac 12}) ).
$$
Further $f(s,d)$ has no zeros in Re $s > 1+ 3 \log \log M/\log M$, and
taking $\sigma_0 = 1+ 3\log \log M/\log M$, we have
$$
{\Cal S}(I_3(d);\Phi) \ll \exp\Big( \pi \frac{(1/2+\epsilon)\log X}{2S}\Big)
M^{-(1-b)} X^{\epsilon}.
$$
\endproclaim

The implied constants above, and elsewhere, may depend upon $\epsilon$, and
the polynomial $P$.
Proposition 2.1 allows us to mollify a little away from $\frac 12$, and
we now turn to the more delicate question of mollifying near $\frac 12$.
Let $\delta_1$ and $\delta_2$ be two complex numbers, and
define $\tau= \frac{\delta_1+\delta_2}2$, and $\delta =
\frac{\delta_1-\delta_2}{2}$.  We put $\kappa=
\max(|\text{Re }\delta|, |\text{Re }\tau|)$
and we suppose below that $\kappa \le \frac 14$.  Let
$$
\xi(s,\chi_{-8d})=\biggl(\frac{8d}{\pi}\biggr)^{\frac{s}{2}-\frac 14}
\Gamma(\tfrac s2+ \tfrac 12) L(s,\chi_{-8d}),
$$
denote the completed $L$-function which satisfies the functional equation
$\xi(s)=\xi(1-s)$.  We shall show how to evaluate
$$
{\Cal S}(\xi(\tfrac 12+\delta_1,\chi_{-8d}) \xi(\tfrac 12+\delta_2,\chi_{-8d})
M(\tfrac 12+\delta_1,d)M(\tfrac 12+\delta_2,d);\Psi),
\tag{2.9}
$$
where $\Psi$ is a smooth function supported on $[1,2]$.
We obtain $(8X/\pi)^{\tau}\Gamma_\delta(\tau){\Cal W}(\delta_1,\Phi)$
by taking $\delta_2=\overline{\delta_1}$,
and taking $\Psi(t)= \Phi(t)t^{-\tau}$.

To evaluate the expression
 (2.9), we first need an ``approximate functional equation''
for $\xi(\frac 12+\delta_1,\chi_{-8d})\xi(\frac 12+\delta_2,\chi_{-8d})$.
Put
$$
\Gamma_\delta(s) = \Gamma\Big(\frac 34 +\frac s2 +\frac{\delta}2\Big)
\Gamma\Big(\frac 34+\frac s2 -\frac {\delta}{2}\Big).
$$
For $\xi >0$ we define
$$
W_{\delta,\tau}(\xi) = \frac{1}{2\pi i} \int_{(c)} \Gamma_\delta(s)
\xi^{-s} \frac{2s}{s^2 -\tau^2} ds, \tag{2.10}
$$
where $c > |\text{Re } \tau|$ is a real number.  Here, and throughout,
we abbreviate $\int_{c-i\infty}^{c+i\infty}$ to $\int_{(c)}$.
In Lemma 3.1 we shall see that $W_{\delta,\tau}(\xi)$ is a smooth
function on $(0,\infty)$, and that it decays exponentially as $\xi \to
\infty$.
For all integers $n\ge 1$, and complex numbers $s$ we
put
$$
r_s(n) = \sum_{ab=n } \Big( \frac ab\Big)^s,
$$
which is plainly an even function of $s$.
Finally, for all integers $d>0$ we define
$$
A_{\delta,\tau} (d) := \sum_{n=1}^{\infty} \frac{r_{\delta}(n)}{\sqrt{n}}
\L{-8d}{n} W_{\delta,\tau}\L{n\pi }{8d}.
\tag{2.11}
$$
We establish in Lemma 3.2 that for fundamental discriminants $-8d (<0)$
we have
$$
\xi(\frac 12+\delta_1,\chi_{-8d})\xi(\frac 12+\delta_2,\chi_{-8d})
= A_{\delta,\tau}(d).
$$
Thus our expression in (2.9)
becomes ${\Cal S}(A_{\delta,\tau}(d) M(\frac 12+\delta_1,d)M(\frac 12+
\delta_2,d);\Psi)$.

Let $\sqrt{2X} \ge Y>1$ be a real parameter to be chosen later and write
$\mu^2(d) =M_Y(d) + R_Y(d)$ where
$$
M_Y(d) = \sum\Sb l^2 |d \\ l\le Y\endSb \mu(l), \qquad \text{and }\qquad
R_Y(d) = \sum\Sb l^2 |d \\ l>Y\endSb \mu(l).
$$
Given a sequence $\{a_n\}_{n=1}^{\infty}$, and a smooth function $F$
supported on $[1,2]$, we define
$$
\Sum_M(a_d;F) = \Sum_{M,X,Y}(a_d;F)
=\frac{1}{X} \sum_{d \text{ odd}} M_Y(d) a_d F\biggl(
\frac dX\biggr),
$$
and
$$
\Sum_R(a_d;F) = \Sum_{R,X,Y}(f_d;F) = \frac{1}{X} \sum_{d \text{ odd}}
\Big|R_Y(d) a_d
F\Big(\frac{d}{X}\Big)\Big|,
$$
so that $\Sum(a_d;F)=\Sum_M(f_d;F) + O(\Sum_R(a_d;F))$.

\proclaim{Proposition 2.2}  Let $\Psi$ be a smooth function supported
on $[1,2]$, with $\Psi(t) \ll 1$.  With notations as above, and supposing
that $M\le \sqrt{X}$, we have
$$
{\Cal S}_R(A_{\delta,\tau}(d)M(\tfrac 12+\delta_1,d)M(\tfrac 12
+\delta_2,d);\Psi) \ll X^{\kappa+\epsilon}
\Big(\frac{1}{Y} +\frac{X^{-\text{Re }\delta_1/2}+X^{-\text{Re }\delta_2/2}}
{Y^{\frac 12}} + X^{-\text{Re }\tau}\Big).
$$
\endproclaim

It remains lastly to evaluate ${\Cal S}_M(A_{\delta,\tau}(d) M(\tfrac 12+
\delta_1,d)M(\frac 12+\delta_2,d);\Psi)$.   We evaluate
more generally ${\Cal S}_M(A_{\delta,\tau}(d) \L{-8d}{l};\Psi)$
for any odd integer $l$.  To state our result, we need a
few more definitions.  For any two complex numbers
$s$ and $w$ we define
$$
Z(s;w) = \zeta(s-2w)\zeta(s) \zeta(s+2w). \tag{2.12}
$$
We write the odd integer $l$ as $l=l_1l_2^2$, where $l_1$ and $l_2$ are
odd, and $l_1$ is square-free.  For a complex number $w$
with $|\text{Re }w|\le \frac 14$, and a complex number
$s$ with $\text{Re }s >\frac 12$ we define
$\eta_w(s;l)= \prod_{p} \eta_{p;w}(s;l)$ where
$\eta_{2;w}(s;l)= (1-2^{-s-2w})(1-2^{-s})(1-w^{-s+2w})$ and for primes
$p \ge 3$ we have
$$
\eta_{p;w}(s;l) =
\cases
\L{p}{p+1} \Big(1-\frac{1}{p^s}\Big) \Big(1+\frac 1p +\frac 1{p^s}
- \frac{p^{2w}+p^{-2w}}{p^{s+1}} +\frac{1}{p^{2s+1}}\Big)
&\text{if } p\nmid l\\
\L{p}{p+1} \Big(1-\frac{1}{p^s}\Big) &\text{if } p|l_1\\
\L{p}{p+1} \Big(1-\frac{1}{p^{2s}}\Big) &\text{otherwise}.
\\
\endcases
\tag{2.13}
$$
Note that $\eta_w(s;l)$ is absolutely convergent in the
range of our definition.

\proclaim{Proposition 2.3}  With notations as above, we may write
$$
\align
{\Cal S}_M\Big(\L{-8d}{l} A_{\delta,\tau}(d);\Psi\Big)&\\
= \frac{2}{3\zeta(2) \sqrt{l_1}}
\sum_{\mu=\pm} &\Big( r_\delta(l_1)\Gamma_\delta(\mu\tau) \L{8X}{l_1
\pi}^{\mu \tau} {\check \Psi}(\mu \tau) Z(1+2\mu\tau;\delta)
\eta_\delta(1+2\mu\tau;l)\\
&+ r_\tau(l_1) \Gamma_\tau(\mu\delta) \L{8X}{l_1\pi}^{\mu\delta} {\check \Psi}
(\mu\delta)
Z(1+2\mu\delta;\tau)\eta_\tau(1+2\mu\delta;l)\Big)
\\
&+ {\Cal R}(l) + O\Big(\frac{|r_\delta(l_1)|
X^{\kappa+\epsilon}}{Yl_1^{\frac 12+
\kappa}} + \frac{|r_\delta(l_1)|X^{\epsilon}}{(Xl_1)^{\frac 14}}
+ \frac{X^{|\text{Re }\delta|+\epsilon}l_1^{\kappa+|\text{Re }\delta|
-\frac 12}}{Y^{1-2\kappa-2|\text{Re }\delta|}}\Big).
\\
\endalign
$$
Here ${\Cal R}(l)$ is a remainder term bounded for each individual $l$
by
$$
|{\Cal R}(l)| \ll  \frac{l^{\frac 12+\epsilon}Y^{1+\epsilon}}{X^{\frac 12-
|\text{Re }\delta|-\epsilon}} \Psi_{(2)} \Psi_{(3)}^{\epsilon},
$$
and bounded on average by
$$
\sum_{l=L}^{2L-1} |{\Cal R}(l)| \ll \frac{L^{1+\epsilon}Y^{1+\epsilon}}
{X^{\frac 12-|\text{Re }\delta|-\epsilon} } \Psi_{(2)}\Psi_{(3)}^{\epsilon}.
$$
\endproclaim

We shall prove Proposition 2.3 in Section 5.  Observe that although
each of the four main terms in Proposition 2.3 has singularities
(for example the first term has poles when $\tau=0$, or when $\tau=\pm
\delta$), their sum is regular.

Plainly Propositions 2.2 and 2.3 can be used to evaluate the
quantity in (2.9).  However, carrying this out is complicated, and
in an effort to keep the exposition simple we shall restrict our values
of $\delta_1$ and $\delta_2$ to those necessary in evaluating $J_1(X;\Phi)$
and $J_2(X;\Phi)$.

\proclaim{Proposition 2.4}  Let $\Phi$ be a non-negative
smooth function on $[1,2]$
with $\Phi(t) \ll 1$, and with $\int_1^2 \Phi(t) dt\gg 1$.
Let $\delta_1$ be a complex number such
that $\frac 14 \ge \text{Re }\delta_1
\ge -\frac{1}{\epsilon \log X}$, and with $|\delta_1|
\ge \frac{\epsilon}{\log X}$.  We take $\delta_2=\overline{\delta_1}$
so that $\tau= \text{Re }\delta_1$, and $\delta= i\text{Im } \delta_1$.
Then with the
mollifier coefficients as in (2.8), and with $M=X^{\frac 12-\epsilon}$
we have that ${\Cal W}(\delta_1,\Phi)$ equals
$$
%\align
%{\Cal W}(\delta_1,\Phi)
%&=
1+ \Big( \frac{1-(8X/\pi)^{-2\tau}}{2\tau \log M}
-\L{8X}{\pi}^{-\tau} \frac{(8X/\pi)^{\delta}-(8X/\pi)^{-\delta}}{2\delta
\log M}\Big)
%\\&\hskip 1in \times
\int_0^b M^{-2\tau(1-x)} \Big| Q^{\prime}(x) + \frac{Q^{\prime \prime}(x)}{
2\delta_1 \log M}\Big|^2 dx
$$
with an error $O(X^{-\tau-\epsilon} \Phi_{(2)}
\Phi_{(3)}^{\epsilon}+M^{-2\tau(1-b)}|\delta_1|^6 \log^5 X)$.
\endproclaim

We emphasize that the
conditions on $|\delta_1|$ and $\text{Re }\delta_1$ were
assumed only to ease our exposition. In fact, the stated result
holds without these restraints.
Armed with these results, we complete the proof of Theorem 1.

\demo{Proof of Theorem 1}  We take $\Phi$ to be a smooth function supported
in $(1,2)$ such that $\Phi(t) \in [0,1]$ for all $t$, $\Phi(t)=1$ for
$t\in (1+\epsilon,2-\epsilon)$, and $|\Phi^{(\nu)}(t)| \ll_{\nu}
(1/\epsilon)^\nu$ (so that $\Phi_{(2)} \ll 1/\epsilon$, and $\Phi_{(3)}
\ll 1/\epsilon^2$).  Our mollifier is chosen as in (2.8), with $M=X^{\frac 12
-\epsilon}$.  Further we take $\sigma_0$ as in Proposition 2.1, and
$S=\pi/(2(1-b)) + 10\epsilon$.  Using Proposition 2.1 in (2.6) we get that
$$
{\Cal N}(X,\Phi) \le \frac{X{\Cal S}(1;\Phi)}{8S \sinh (\frac{\pi R}{2S})}
(J_1(X;\Phi) + J_2(X;\Phi)) +o(X),
$$
where $J_1$ and $J_2$ are given in (2.7a,b).

Applying Proposition 2.4 we get that for real numbers $u$ and $v$
with $\max(|u|,|v|) \ge \epsilon$ and $\frac 14\log X \ge u\ge -1/\epsilon$
$$
{\Cal W}\Big(\frac{u+iv}{\log X},\Phi\Big)
= {\Cal V}(u,v)
+ O\Big(
%X^{-\epsilon} e^{-u} +
 M^{-2u(1-b)/\log X} \frac{(1+|u|+|v|)^6}{\log X}
\Big),
$$
where
$$
{\Cal V}(u,v):= 1+ \frac{e^{-u}\log X}{2\log M} \Big( \frac{\sinh u}{u}
-\frac{\sin v}{v} \Big) \int_0^b M^{-2u(1-x)/\log X} \Big| Q^{\prime}(x)
+ \frac{Q^{\prime\prime}(x) \log X}{2 (x+iy)\log M} \Big|^2 dx.
$$
Plainly ${\Cal V}(u,v) \ge 1$ always, and so
we deduce that
$$
J_1(X;\Phi) = \int_0^S \cos \L{\pi t}{2S}
\log {\Cal V}(-R,t) dt + O\Big( \frac{1}{\log X}\Big).
$$
Further, using the above together with Proposition 2.1, and
keeping in mind our choice for $S$,  we obtain
that
$$
J_2(X; \Phi) = \int_0^{\infty} \sinh \L{\pi u}{2S}
\log {\Cal V}(u-R, S) du + o(1).
$$
We conclude that
$$
\align
{\Cal N}(X,\Phi) &\le \frac{X{\Cal S}(1;\Phi)}{8S \sinh \L{\pi R}{2S} }
\Big( \int_0^S \cos \L{\pi t}{2S}
\log {\Cal V}(-R,t) dt \\
&\hskip 1.5 in +  \int_0^{\infty} \sinh \L{\pi u}{2S}
\log {\Cal V}(u-R, S) du\Big) + o(X).
\\
\endalign
$$

We now take $R=6.8$, $b=0.64$, and $P(x) =3 (x/b)^2 -2 (x/b)^3$.
Then a computer calculation showed that
${\Cal N}(X,\Phi) \le 0.79 X {\Cal S}(1;\Phi) + o(X)$.
Taking $X=x/2$, $x/4$, $\ldots$, we obtain Theorem 1.

\enddemo

We end this section by reflecting on some features of
the method used to prove Theorem 1.
Our overall strategy was to estimate on average the number of zeros
(weighted suitably) of the mollified $L$-function in a small box ${\Cal B}$
as in Lemma 2.1. If we use the usual argument principle to estimate
the zeros in ${\Cal B}$, then we face the problem of
trying to understand the argument of $f(s,d)$ on the horizontal
sides of ${\Cal B}$.  This appears to be difficult
because the argument of $L(s,\chi_{-8d})$ is intimately related
to the location of its zeros.  Selberg's argument principle (Lemma 2.1)
allows us to circumvent this by introducing the kernel
$\sin(\pi (s-W_0+iH)/(2iH))$ which is real on the left vertical
edge of ${\Cal B}$, and purely imaginary on the horizontal edges of
${\Cal B}$.  This enables us to deal only with $\log |f(s,d)|$
(a quantity well suited for estimating from above) on these three
sides of ${\Cal B}$, while on the left vertical edge of ${\Cal B}$
we are in the region of absolute convergence of $L(s,\chi_{-8d})$
so that $\log f(s,d)$ is relatively easy to understand on this line.

The chief drawback with Selberg's lemma is the exponential
growth of the the kernel $\sin(\pi(s-W_0+iH)/(2iH))$ on the horizontal
sides of ${\Cal B}$.  To offset this it is necessary that $\log |f(s,d)|$
be very small on the horizontal sides of ${\Cal B}$
(at least on average over $d$).  This motivates our choice (see (2.8))
of the mollifier coefficients $\lam(n)$:  this choice
guarantees that the Dirichlet series coefficients of
$f(s,d)$ vanish for $2\le n \le M^{1-b}$ so that we would expect
$f(s,d)$ to be close to $1$ on average (as confirmed by Proposition
2.1).  Since the growth of Selberg's kernel is determined by the
height of the box $H$, and the decay of $\log |f(s,d)|$ is
controlled by how long a mollifier we can take, we see that
there is a natural limitation on how small a box we can take
in terms of how long a mollifier we can allow.

In this way we reduce the problem of estimating the
weighted average of zeros in ${\Cal B}$ to evaluating certain
mollified mean values, and that is accomplished by extending the
ideas in [11].  There are two features of this approach
which are a little dissatisfying.  Firstly the choice of mollifier
coefficients is made in an {\sl ad hoc} way through some
numerical experimentation.  This is in contrast with the
classical situation of mollifying at a point where the optimal
mollifier coefficients emerge as minimizers of a certain quadratic
form while keeping a linear form fixed.  The situation here
is less clear because the final answer depends on a complicated integral
over the sides of ${\Cal B}$ of the mollified moments, and also
because the initial mollifier coefficients are no longer free, as explained
above.  We have not understood this optimization problem fully, and
it is quite possible that a better choice of mollifier exists.

Secondly, the proof of Theorem 1 relied crucially
upon knowing that our weighted average of zeros is less than $1$.
Since this emerged only after an involved calculation we now indicate
why it is reasonable to expect this average to be small.
More precisely note that in the proof of Theorem 1 we bounded
$$
\frac{1}{2\sinh \L{\pi R}{2S}} \sum\Sb \beta+i\gamma \in {\Cal B} \\
f(\beta+i\gamma,d)=0\endSb \cos \L{\pi \gamma \log X}{2S}
\sinh\L{\pi(R+\log X(\beta -1/2))}{2S}. \tag{2.14}
$$
We showed that on average over $d$ this quantity is bounded by $0.79$,
while if $L(s,\chi_{-8d})$ had a real zero this quantity exceeds $1$;
thus producing many $L(s,\chi_{-8d})$ having no real zeros.  We now
restrict our attention to the zeros in (2.14) arising from $L(s,\chi_{-8d})$
term, and calculate (conjecturally) their contribution.  We
suspect that the contribution from zeros of the
mollifier  to (2.14) is negligible on average; at any rate
(2.14) is at least as large as the contribution from zeros
of $L(s,\chi_{-8d})$, and so it is necessary that this be small.
If we assume the Generalized Riemann Hypothesis then the
zeros of $L(s,\chi_{-8d})$ in ${\Cal B}$
contribute to (2.14) the amount
$$
\frac{1}{2} \sum\Sb |\gamma\log X|\le S \\ L(\tfrac 12+i\gamma,\chi_{-8d})=0
\endSb \cos \L{\pi \gamma\log X}{2S}. \tag{2.15}
$$
The distribution of low lying zeros in families of $L$-functions
has been studied extensively by Katz and Sarnak [7], and the
conjectures they formulate there enable one to calculate
sums like (2.15) on average.  Our family of $L$-functions
is expected to have a symplectic symmetry, whose $1$-level density
function is conjecturally $1-\sin(2\pi x)/(2\pi x)$ (see pages 405-409 of
[7]).  Note that this density vanishes to order $2$ at $0$, indicating that
the zeros of $L(s,\chi_{-8d})$ tend to repel the point $1/2$.
This philosophy predicts that the
average value of (2.15) is
$$
\int_{0}^{S/(2\pi)} \cos \L{\pi^2 x}{S} \Big(1-\frac{\sin(2\pi x)}{2\pi x}
\Big) dx.
$$
For the choice of $S$ in Theorem 1 (namely $S = \pi/(0.72)$) the above
evaluates to $0.1827\ldots$.  Thus conjecturally there
are very few zeros in our box, and this suggests an explanation
for why the method works.

We may ask if results similar to Theorem 1 hold for other
families of $L$-functions.  Our remarks above
indicate that perhaps the method would succeed in other
families with a repulsion phenomenon at $1/2$.  One example
of these is the family of modular forms (say, of large weight) and
odd sign of the functional equation, where there is always a zero
at $1/2$ but the next zero is repelled.  We hope to return to
these questions later.

\head 3. Preliminaries \endhead

\subhead 3.1 The approximate functional equation \endsubhead

\noindent

\proclaim{Lemma 3.1}  For $\xi \in (0,\infty)$, $W_{\delta,\tau}(\xi)$
is a smooth complex-valued function.  For $\xi$ near $0$ we
have the asymptotic
$$
W_{\delta,\tau}(\xi) = \Gamma_\delta(\tau)\xi^{-\tau} +\Gamma_\delta(-\tau)
\xi^{\tau} + O(\xi^{1-\epsilon}).
$$
For large $\xi$ and any integer $\nu$ we have the estimate
$$
W_{\delta, \tau}^{(\nu)} (\xi) \ll_{\nu} \xi^{2\nu+6} e^{-2\xi}
\ll_{\nu} e^{-\xi}.
$$
\endproclaim

\demo{Proof}  By moving the line of integration in (2.10)
to Re $s = -1+\epsilon$
we see immediately the asymptotic claimed for small $\xi$.  Plainly
the $\nu$-th derivative of $W_{\delta,\tau}$ is given by
the convergent integral
$$
\frac{(-1)^{\nu}}{2\pi i} \int_{(c)} \Gamma_\delta (s) s(s+1)\cdots
(s+\nu-1) \xi^{-s} \frac{2s}{s^2 -\tau^2} ds
$$
for any $c>|\text{Re }\tau|$.  Thus $W_{\delta,\tau}(\xi)$ is smooth.
To prove the last estimate of the lemma we may suppose that
$\xi > \nu +4$.  Since $|\Gamma(x+iy)|\le \Gamma(x)$ for $x\ge 1$, and
$s\Gamma(s)=\Gamma(s+1)$, we obtain that the integral above giving
$W_{\delta,\tau}^{(\nu)}(\xi)$ is bounded by
$$
\ll_{\nu} |\Gamma(c/2+\nu+3)|^2 \xi^{-c} \int_{(c)}
\frac{|ds|}{|s^2-\tau^2|} \ll_{\nu} \Gamma(c/2 +\nu+3)^2 \frac{\xi^{-c}}
{c-|\text{Re }\tau|}.
$$
By Stirling's formula this is
$$
\ll_{\nu} \Big( \frac{c+2\nu+6}{2e}\Big)^{c+2\nu+6} \frac{\xi^{-c}}
{c-|\text{Re }\tau|},
$$
and taking $c= 2\xi -2\nu-6 (\ge 2)$ we get the lemma.

\enddemo

\proclaim{Lemma 3.2}  Recall that
$\kappa = \max(|\text{Re }\delta|,|\text{Re }\tau|)\le \frac 14$.
For fundamental discriminants $-8d (<0)$ we have
$$
\xi(\tfrac 12+\delta_1,\chi_{-8d})\xi(\tfrac 12+\delta_2,\chi_{-8d}
) = A_{\delta,\tau}(d).
$$
\endproclaim
\demo{Proof}  Consider for some $3/2 - |\text{Re }\delta| > c > 1/2+
|\text{Re }\delta|$
$$
\frac{1}{2\pi i}\int_{(c)} \xi(\tfrac 12+\delta+s,\chi_{-8d})
\xi(\tfrac 12-\delta+s,\chi_{-8d}) \frac{2s}{s^2 -\tau^2} ds.
$$
Expanding $L(\tfrac 12+\delta+s,\chi_{-8d}) L(\tfrac 12-\delta+s,\chi_{-8d})$
into its Dirichlet series
$\sum_{n=1}^{\infty}\frac{ r_\delta(n)}{n^{\frac 12+s}} \L{-8d}{n}$,
and integrating term by term,
we get that this equals $A_{\delta,\tau}(d)$.
Now move the path of integration to the line Re$(s)=-c$.  We encounter
poles at $s=\tau$, $-\tau$, and the residues here give
$\xi(\frac 12+\delta+\tau,\chi_{-8d})\xi(\frac 12-\delta+\tau,\chi_{-8d})
+\xi(\frac 12+\delta-\tau,\chi_{-8d})\xi(\frac 12-\delta-\tau,\chi_{-8d})
= 2\xi(\frac 12+\delta_1,\chi_{-8d})\xi(\frac 12+\delta_2,\chi_{-8d})$,
upon using the functional equation.
In the remaining integral on the $-c$ line,
we let $s \to -s$ and use the functional equation.  Then it
evaluates to $-A_{\delta,\tau}(d)$, which completes our proof.

\enddemo

\subhead 3.2 On Gauss-type sums \endsubhead

\noindent Let $n$ be an odd integer.  We define for all integers $k$
$$
G_k(n) = \biggl(\frac{1-i}{2} +\fracwithdelims() {-1}{n} \frac{1+i}{2}
\biggr) \sum_{a\pmod n}\fracwithdelims() an e\biggl(\frac{ak}{n}\biggr),
$$
and put
$$
\tau_k(n) = \sum_{a\pmod n} \L{a}{n} e\L{ak}{n} = \biggl(\frac{1+i}{2}
+\L{-1}{n} \frac{1-i}{2}\biggr) G_k(n).
$$
If $n$ is square-free then $\L{\cdot}{n}$ is a primitive character with
conductor $n$.  Here it is easy to see that $G_k(n) = \L{k}{n} \sqrt{n}$.
For our later work, we require knowledge of $G_k(n)$ for all odd $n$.  This
is contained in the next Lemma which is Lemma 2.3 of [11].

\proclaim{Lemma 3.3}
\item{(i)} (Multiplicativity)  Suppose $m$ and $n$ are coprime odd integers.

\noindent Then $G_k(mn) =G_k(m) G_k(n)$.
\item{(ii)} Suppose $p^{\alpha}$ is the largest power of $p$ dividing $k$.
(If $k =0$ then set $\alpha =\infty$.)  Then for $\beta \ge 1$
$$
G_k(p^{\beta})=
\cases
0&\text{if   } \beta \le \alpha \text{  is odd},\\
\phi(p^\beta) & \text{if   } \beta \le \alpha \text{   is even},\\
-p^{\alpha} &\text{if   } \beta = \alpha +1 \text{   is even},\\
\fracwithdelims(){kp^{-\alpha}}{p} p^{\alpha} \sqrt p &\text{if   } \beta
=\alpha +1 \text{   is odd}\\
0&\text{if   } \beta \ge \alpha +2.\\
\endcases
$$
\endproclaim

\subhead 3.3 Lemmas for estimating character sums \endsubhead

\noindent We collect here two lemmas that will be very useful in
bounding the character sums that arise below.  These are
consequences of a recent large sieve result for real characters
due to D. R. Heath-Brown [4] (see Lemmas 2.4 and 2.5 of [11]).

\proclaim{Lemma 3.4}  Let $N$ and $Q$ be positive integers and let
$a_1$, $\ldots$, $a_N$ be arbitrary complex numbers.  Let $S(Q)$ denote
the set of real, primitive characters $\chi$ with conductor $\le Q$.  Then
$$
\sum_{\chi \in S(Q)} \biggl| \sum_{n\le N} a_n\chi(n)  \biggr|^2
\ll_{\epsilon}  (QN)^{\epsilon} (Q+N)\sum_{n_1n_2=\square} |a_{n_1}a_{n_2}|,
$$
for any $\epsilon >0$.  Let $M$ be any positive integer, and
for each $|m|\le M$ write $4m= m_1 m_2^2$ where $m_1$ is a
fundamental discriminant, and $m_2$ is positive.  Suppose
the sequence $a_n$ satisfies $|a_n|\ll n^{\epsilon}$.  Then
$$
\sum_{|m| \le M} \frac{1}{m_2} \
 \biggl| \sum_{n\le N} a_n \L{m}{n} \biggr|^2 \ll (MN)^{\epsilon} N(M+N).
$$
\endproclaim

\proclaim{Lemma 3.5}  Let $S(Q)$ be as in Lemma 3.4, and suppose
$\sigma+it$ is a complex number with $\sigma\ge \frac 12$.
Then
$$
\sum_{\chi \in S(Q)}  |L(\sigma+it,\chi)|^4
\ll Q^{1+\epsilon}(1+|t|)^{1+\epsilon}, \ \ \
\text{and} \ \ \
\sum_{\chi \in S(Q)} |L(\sigma+it,\chi)|^2 \ll Q^{1+\epsilon}
(1+|t|)^{\frac 12+\epsilon}.
$$
\endproclaim

\subhead 3.4 Poisson summation \endsubhead

\noindent For a Schwarz class function $F$ we define
$$
{\tilde F}(\xi) = \frac{1+i}{2}{\hat F}(\xi) + \frac{1-i}{2}{\hat F}(-\xi)
= \int_{-\infty}^{\infty} (\cos (2\pi \xi x) + \sin (2\pi \xi x)) F(x) dx.
$$
We quote the following version of Poisson summation (see Lemma 2.6 of [11]):

\proclaim{Lemma 3.6}  Let $F$ be a smooth function
supported in $(1,2)$.  For any odd integer $n$,
$$
\Sum_M\biggl(\L{d}{n};F\biggr) =  \frac{1}{2n} \L {2}{n}
\sum\Sb \alpha \le Y \\ (\alpha, 2n)=1\endSb \frac{\mu(\alpha)}{\alpha^2}
\sum_k (-1)^k G_k(n) {\tilde F} \L {kX}{2\alpha^2 n}.
$$
\endproclaim

\head 4. Proofs of Propositions 2.1 and 2.2 \endhead

\noindent We first record two applications of Lemma 3.4 which
will be useful in the proofs of these Propositions.
Write  $\lambda_2(n) =\sum_{ab=n, a,b\le M}
\lambda(a)\lambda(b)$.  Note that $|\lambda_2(n)|\ll n^{\epsilon}$ and
that $M(s,d)^2 =\sum_{n\le M^2} \lambda_2(n) n^{-s}
\L{-8d}{n}$.  By Lemma 3.4 we see that
for $N \le M^2 (\ll X^2)$ we
have
$$
\align
\sum_{X\le d\le 2X} \mu^2(2d)
\Big|\sum_{N\le n\le 2N} \frac{\lam_2(n)}{n^s} \L{-8d}{n}\Big|^2
&\ll X^{\epsilon} (X+N) \sum\Sb N\le n_1,n_2 \le 2N \\ n_1n_2=\square
\endSb \frac{|\lam_2(n_1)\lam_2(n_2)|}{(n_1 n_2)^{\text{Re }s}}\\
&\ll X^{\epsilon} (X+N) N^{1-2\text{Re }s}
\sum\Sb N\le n_1,n_2 \le 2N \\ n_1n_2=\square
\endSb \frac{1}{\sqrt{n_1n_2}}\\
&\ll  X^{\epsilon} (X+N) N^{1-2 \text{Re }s} \sum_{a\le M^2} \frac{d(a^2)}{a}
\\
&\ll  X^{\epsilon} (X+N) N^{1-2\text{Re }s}.\\
\endalign
$$
From this we conclude that
$$
\sum_{X\le d\le 2X} \mu^2(2d) |M(s,d)|^4
\ll X^{\epsilon} (X+ XM^{2(1-2\text{Re }s)} +M^{4(1-\text{Re }s)}).
\tag{4.1}
$$
In a similar manner we see that
if $l$ is any odd integer $\le \sqrt{2X}$ then
$$
\sum_{X/l^2 \le m\le 2X/l^2} \mu^2(2m) |M(s,l^2m)|^4
\ll X^{\epsilon} \Big(\frac{X}{l^2} + \frac{X}{l^2} M^{2(1-2\text{Re }s)}
+ M^{4(1-\text{Re }s)}\Big). \tag{4.2}
$$

\subhead 4.1 Proof of Proposition 2.1 \endsubhead

\noindent Since $\Phi$ is a non-negative smooth function supported on $[1,2]$
such that  $\Phi(t) \ll 1$, and $\int_1^2 \Phi(t) dt \gg 1$
we see that ${\Cal S}(1;\Phi) \gg X^{-1} \sum_{X\le d \le 2X} \mu^2(2d)
\gg 1$.  We write $B(s,d) = L(s,\chi_{-8d}) M(s,d) -1$ so that
$$
{\Cal W}(\delta_1,\Phi) = 1+ O({\Cal S}(B(\tfrac 12+\delta_1,d);\Phi)
+ {\Cal S}(|B(\tfrac 12+\delta_1,d)|^2;\Phi)). \tag{4.3}
$$
To estimate the unknown terms above, we consider
$$
\frac{1}{2\pi i} \int_{(c)} \Gamma(s+1) B(\tfrac 12+\delta_1+s,d) X^s
\frac{ds}{s},
$$
for any real number $c> \frac 12 -\text{Re }\delta_1$.  We move
the line of integration to the line Re $s = -\text{Re }\delta_1$.  The
pole at $s=0$ contributes $B(\frac 12+\delta_1;d)$ and so
we conclude that $B(\tfrac 12+\delta_1,d)$ equals
$$
\frac{1}{2\pi i} \int_{(c)}
\Gamma(s+1) B(\tfrac 12+\delta_1+s,d) X^s \frac{ds}{s}
- \frac{1}{2\pi i}\int_{(-\text{Re }\delta_1)}
\Gamma(s+1) B(\tfrac 12+\delta_1+s,d) X^s \frac{ds}{s}. \tag{4.4}
$$
Write the expression in (4.4) as $T_1(\frac 12+\delta_1,d)-
T_2(\frac 12+\delta_1,d)$, say.

We first consider the contributions of the $T_2(d)$ terms
to the unknown quantities in (4.3).  We shall prove
that
$$
{\Cal S}(|T_2(\tfrac 12+\delta_1, d)|^2; \Phi) \ll
X^{-2\text{Re }\delta_1 +\epsilon},
\,\,\,\text{and } \,\,\,{\Cal S}(|T_2(\tfrac 12+\delta_1,d)|;\Phi)
\ll X^{-\text{Re }\delta_1 +\epsilon}. \tag{4.5}
$$
Plainly the second estimate above follows from the first and Cauchy's
inequality.  To see the first estimate observe that by
Cauchy's inequality
$$
|T_2(\tfrac 12+\delta_1,d)|^2 \ll X^{-2\text{Re }\delta_1}
\Big( \int_{(-\text{Re }\delta_1)} |\Gamma(s+1)B(\tfrac 12+\delta_1+s,d)^2 ds|
\Big)\Big(\int_{(-\text{Re }\delta_1)} |\Gamma(s+1)| \frac{|ds|}{|s|^2}\Big),
$$
and in view of the rapid decay of $|\Gamma(s+1)|$ as $|\text{Im s}| \to
\infty$, we deduce that
$$
|T_2(\tfrac 12+\delta_1,d)|^2 \ll X^{-2\text{Re }\delta_1}
\Big( 1+ \int_{(-\text{Re }\delta_1)} |\Gamma(s+1)| |L(\tfrac 12+\delta_1+s,
\chi_{-8d})
M(\tfrac 12+\delta_1+s,d)|^2 |ds|\Big).
$$
Averaging this over the appropriate $d$, with another
application of Cauchy's inequality
we obtain that ${\Cal S}(|T_2(\frac 12+\delta_1,d)|^2;\Phi)$ is bounded
by
$$
 X^{-2\text{Re }\delta_1} \Big( 1+
 \int_{(-\text{Re }\delta_1)} |\Gamma(s+1)| {\Cal S}(|L(\tfrac 12+\delta_1+s,
\chi_{-8d})|^4;\Phi)^{\frac 12}
{\Cal S}(|M(\tfrac 12+\delta_1+s,d)|^4;\Phi)^{\frac 12} |ds| \Big),
$$
and (4.5) follows upon using Lemma 3.5 and (4.1) above (keeping
in mind that $M\le \sqrt{X}$).

It remains now to consider the $T_1$ contribution.  In the region Re $s >1$
we may write
$$
B(s,d) = \sum_{n=1}^{\infty} \frac{b(n)}{n^s} \L{-8d}{n}.
$$
From the shape of our mollifier we see that
$b(n)=0$ for all $n\le M^{1-b}$, $b(n)=0$ for all
square values $n \le M^{2(1-b)}$ (because $b(m^2) = \sum_{d|m^2}
\lam(d) =\sum_{d|m}\lam(d)=b(m)$, since $\lam$ is supported
on square-free numbers), and lastly $|b(n)|\ll d(n)\ll n^{\epsilon}$
for all $n$.  We write
$$
\align
T_1(\tfrac 12+\delta_1,d) = \frac{1}{2\pi i} \int_{(c)}&\Gamma(s+1)
\sum_{M^{1-b} \le n\le X\log^2 X} \frac{b(n)}{n^{\frac 12+\delta_1+s}}
\L{-8d}{n} X^s \frac{ds}{s}\\
&+ \sum_{n> X\log^2 X} \frac{b(n)}{n^{\frac 12+\delta_1}} \L{-8d}{n}
\Big(\frac{1}{2\pi i}\int_{(c)} \Gamma(s+1) \L{X}{n}^s
\frac{ds}{s}\Big).\\
\endalign
$$
In the second term above we take $c=n/(10X)$ and use Stirling's formula
to get that $\frac{1}{2\pi i} \int_{(c)} \Gamma(s+1) \L{X}{n}^s \frac{ds}{s}
\ll \exp(-n/(20X))$.  Hence the second term above
contributes $\ll X^{-5}$ say.  In  the first term above we move the
line of integration to Re $s=\frac 1{\log X}$.  Thus
$$
T_1(\tfrac 12+\delta_1,d) =
 \frac{1}{2\pi i} \int_{(\frac{1}{\log X})}
\Gamma(s+1)
\sum_{M^{1-b} \le n\le X\log^2 X} \frac{b(n)}{n^{\frac 12+\delta_1+s}}
\L{-8d}{n} X^s \frac{ds}{s} +(X^{-5}).
\tag{4.6}
$$

By Cauchy's inequality we get that
$$
\align
|T_1(\tfrac 12+\delta_1,d)|^2 &\ll X^{-10} +
\Big(\int_{(\frac 1{\log X})} |\Gamma(s+1)|\Big| \sum_{M^{1-b} \le n\le
X\log^2 X} \frac{b(n)}{n^{\frac 12+\delta_1+s}} \L{-8d}{n}\Big|^2 |ds|
\Big)\\
&\hskip 1 in \times
 \Big(\int_{(\frac 1{\log X})} |\Gamma(s+1)| \frac{|ds|}{|s^2|}\Big)\\
&\ll X^{-10} + X^{\epsilon}
\int_{(\frac 1{\log X})} |\Gamma(s+1)|\Big| \sum_{M^{1-b} \le n\le
X\log^2 X} \frac{b(n)}{n^{\frac 12+\delta_1+s}} \L{-8d}{n}\Big|^2 |ds|.
\\
\endalign
$$
Splitting the sum over $n$ into dyadic blocks and using Lemma
3.4 we conclude that
$$
{\Cal S}(|T_1(\tfrac 12+\delta_1,d)|^2;\Phi)
\ll  M^{-2 \text{Re }\delta_1(1-b)}X^{\epsilon},
$$
which when combined with (4.5) gives that
$$
{\Cal S}(|B(\tfrac 12+\delta_1,d)|^2;\Phi) \ll M^{-2\text{Re }\delta_1(1-b)}
X^{\epsilon}. \tag{4.7}
$$

We now show how to bound
${\Cal S}(T_1(\tfrac 12+\delta_1,d);\Phi)$.  By (4.6) we see that
$$
{\Cal S}(T_1(\tfrac 12+\delta_1,d);\Phi) \ll
X^{-5} + X^{\epsilon} \sum_{M^{1-b}\le n\le X\log^2 X} \frac{|b(n)|}
{n^{\frac 12+\text{Re }\delta_1}}
\Big|{\Cal S}\Big( \L{-8d}{n};\Phi\Big)\Big|. \tag{4.8}
$$
For each odd integer $n$ let $\psi_n$ denote the character $\psi_n(m)
= \L{m}{n}$ whose conductor is at most $n$.  Note that $\psi_n$ is
non-trivial unless $n$ is a square.  Observe that
for any sequence of numbers $a_n \ll n^{\epsilon}$, and any smooth
function $g$ with $g(0)=0$ and $g(x)$ decaying rapidly as $x\to \infty$,
we have the Mellin transform identity
$$
\sum_{n=1}^{\infty} a_n g(n) = \frac{1}{2\pi i} \int_{(c)} \sum_{n=1}^{\infty}
\frac{a_n}{n^w} \Big( \int_0^\infty g(t) t^{w-1} dt\Big) dw, \tag{4.9}
$$
where $c>1$.  Hence we obtain that
for any odd integer $n$
$$
\align
{\Cal S}(\psi_n(-8d);\Phi) &= \frac{\psi_n(-8)}{2\pi i}
\int_{(c)} \sum_{d=1}^{\infty} \frac{\mu^2 (2d) \psi_n(d) }{d^w}
X^{w-1} {\check \Phi}(w-1) dw \\
&=  \frac{\psi_n(-8)}{2\pi i}
\int_{(c)} \frac{L(w,\psi_n)}{L(2w,\psi_n)} (1+\psi_n(2)/2^w)^{-1}
X^{w-1} {\check \Phi}
(w-1)dw,\\
\endalign
$$
where $L(w,\psi_n)= \sum_{d=1}^{\infty} \psi_n(d)/d^w$ is the
usual Dirichlet $L$-function.  We move the line of integration
above to the line Re $w = \frac 12+ \frac{1}{\log X}$.  We
encounter a pole at $w=1$ if and only if $n$ is a square
(in which case $L(w,\psi_n)$ is essentially $\zeta(w)$) and the
residue of this pole is $\ll 1$.  Thus
we conclude that
$$
|{\Cal S}(\psi_n(-8d);\Phi)| \ll \delta(n=\square) + X^{-\frac 12 +
\epsilon}
\int_{(\frac 12+\frac
{1}{\log X})} |L(w,\psi_n)| |{\check \Phi}(w-1)| |dw|,
$$
where $\delta(n=\square)$ is $1$ if $n$ is a square, and $0$ otherwise.
Since $b(n)=0$ for all squares $\le M^{2(1-b)}$ we
find that
$$
\align
\sum_{M^{1-b}\le n\le X\log^2 X} &\frac{|b(n)|}
{n^{\frac 12+\text{Re }\delta_1}}
\Big|{\Cal S}\Big( \L{-8d}{n};\Phi\Big)\Big|
\ll X^{\epsilon}  M^{-2\text{Re }\delta_1(1-b)}\\
& +
X^{-\frac 12+\epsilon} \int_{(\frac 12+\frac
{1}{\log X})} \sum_{M^{1-b}\le n\le X\log^2 X} \frac{1}
{n^{\frac 12+\text{Re }\delta_1}}  |L(w,\psi_n)| |{\check\Phi}(w-1)||dw|.
\\
\endalign
$$
An easy application of Lemma 3.5 gives that
$$
\sum_{N\le n\le 2N} |L(w,\psi_n)| \ll N^{1+\epsilon} (1+|w|)^{\frac 14
+\epsilon},
$$
and using this above, together with (2.12c) (taking $\nu=2$ there),
we get that
$$
\align
\sum_{M^{1-b}\le n\le X\log^2 X} \frac{|b(n)|}
{n^{\frac 12+\text{Re }\delta_1}}
&\Big|{\Cal S}\Big( \L{-8d}{n};\Phi\Big)\Big|
\\
%&\ll
% X^{\epsilon} \Big( M^{-2\text{Re }\delta_1(1-b)} + X^{-\text{Re }\delta_1}
%+ M^{(\frac 12-\text{Re }\delta_1)(1-b)}X^{-\frac 12}\Big)
%\int_{(\frac 12+\frac{1}{\log X})} (1+|w|)^{\frac 14+\epsilon}
%|{\check\Phi}(w-1)| |dw|
%\\
&\ll
X^{\epsilon} \Phi_{(2)}
\Big( M^{-2\text{Re }\delta_1(1-b)} + X^{-\text{Re }\delta_1}
+ M^{(\frac 12-\text{Re }\delta_1)(1-b)}X^{-\frac 12}\Big).\\
\endalign
$$
Using this in (4.8), and combining with (4.5) we
deduce that (since $M \le \sqrt{X}$)
$$
{\Cal S}(B(\frac 12+\delta_1,d);\Phi)
\ll X^{\epsilon} \Phi_{(2)}
\Big( M^{-2\text{Re }\delta_1(1-b)}
+ M^{(\frac 12-\text{Re }\delta_1)(1-b)}X^{-\frac 12}\Big). \tag{4.10}
$$
Using (4.10) and (4.7) in (4.3), we deduce the
first statement of the Proposition.

To see the second assertion, note that
$$
f(s,d) = 1+B(s,d) = 1+ O\Big( \sum_{n\ge M^{1-b}} \frac{d(n)}{n^{\text{Re }s}}
\Big),
$$
from which it follows easily that $f(s,d)$ has no zeros to the
right of $1+3\log \log M/\log M$. Further
for $s$ in this region $\log f(s,d) = B(s,d) + O(|B(s,d)|^2)$,
and so, with $\sigma_0$ as in the Proposition, we
have
$$
{\Cal S}(I_3(d);\Phi) \ll \exp\Big(\pi \frac{(\frac 12+\epsilon)\log X}{2S}
\Big) \Big(  |{\Cal S}(B(s,d);\Phi)| +{\Cal S}(|B(s,d)|^2;\Phi)\Big).
$$
Thus the second assertion also follows from (4.7) and (4.10).

\subhead 4.2 Proof of Proposition 2.2 \endsubhead

\noindent Observe that $R_Y(d)=0$ unless $d=l^2m$ where $m$
is squarefree and $l>Y$.  Further, note that $|R_Y(d)| \le \sum_{k|d} 1
\ll d^{\epsilon}$.  Hence
$$
\align
\Sum_R(A_{\delta,\tau}(d)&M(\tfrac 12+\delta_1,d)M(\tfrac 12+\delta_1,d)
;\Psi)\\
&\ll X^{-1+\epsilon} \sum\Sb Y < l \\ (l,2)=1\endSb
\ \ \sumflat_{X/l^2 \le m\le 2 X/l^2}
|A_{\delta,\tau}(l^2m) M(\tfrac 12+\delta_1,l^2m)M(\tfrac 12+\delta_2,l^2m)|,
\tag{4.11}\\
\endalign
$$
where the $\flat$ on the sum over $m$ indicates that $m$
is odd and squarefree.  By two applications of
Cauchy's inequality the sum over $m$ above is
$$
\ll \biggl(\sumflat_{ m}|M(\tfrac 12+\delta_1,l^2m)|^{4}
\biggr)^{\frac14}
 \biggl(\sumflat_{m }|M(\tfrac 12+\delta_2,l^2m)|^{4}
\biggr)^{\frac14}
\biggl(\sumflat_{m} |A_{\delta,\tau}(l^2m)|^2
\biggr)^{\frac 12}. \tag{4.12}
$$

Now observe that for any $c>\frac 12 +|\text{Re }\delta|$
$$
\align
A_{\delta,\tau}(l^2 m)
%&=\sum_{n=1}^{\infty} \frac{d_j(n)}{\sqrt{n} } \L
%{8l^2 m}{n} \om_j \biggl(\frac{\pi n}{8 l^2 m}\biggr)\\
&=\frac {1}{2\pi i} \int_{(c)}
\Gamma_\delta(s)
\L{8l^2 m}{\pi}^{s} \frac{2s}{s^2 -\tau^2}
\sum_{n=1}^{\infty} \frac{r_\delta(n)}{n^{s+\frac 12}}\L {-8l^2 m}{n}
ds. \tag{4.13}\\
\endalign
$$
Plainly
$$
\sum_{n=1}^{\infty} \frac{r_\delta(n)}{n^{s+\frac 12}}\L{-8l^2m}{n} =
L(\tfrac 12+s+\delta,\chi_{-8m})L(\tfrac 12 + s-\delta,\chi_{-8m})
 {\Cal E}(s,l)
\tag{4.14}
$$
where
$$
{\Cal E}(s,l)= \prod_{p|l} \biggl(1-\frac{1}{p^{s+\frac 12+\delta}}
\L{-8m}{p}\biggr)\biggl(1-\frac{1}{p^{s+\frac 12-\delta}}\L{-8m}{p}\biggr).
$$
Since $\chi_{-8m}$ is non-principal,
it follows that the left side of (4.14) is analytic for all $s$.

Hence we may move the line of integration
in (4.13) to the line from $\kappa+ 1/\log X-i\infty$ to
$\kappa + 1/\log X+i\infty$.  We encounter
no poles, and so $A_{\delta,\tau}(l^2m)$ is
given by the integral on this new line.  Since
$|{\Cal E}(s,l)|\le \prod_{p|l} (1+1/\sqrt{p})^2 \ll l^{\epsilon}
\ll X^{\epsilon}$,  $2s/(s^2 -\tau^2) \ll X^\epsilon$, and
$|\Gamma_\delta(s)|$ decays exponentially for large $|\text{Im }s|$, we
obtain by Cauchy's inequality that
$$
|A_{\delta,\tau}(l^2m)|^2 \ll X^{2\kappa+\epsilon}
\int_{(\kappa+\frac 1{\log X})}
|\Gamma_\delta(s)| |L(\tfrac 12+s+\delta, \chi_{-8m})L(\tfrac 12+s-\delta,
\chi_{-8m})|^{2}
|ds|.
$$
Summing this over $m$ and using Lemma 3.5, we obtain that
$$
\sumflat_{X/l^2 \le m\le 2X/l^2}
|A_{\delta,\tau}(l^2 m)|^2
\ll \frac{X^{1+2\kappa+\epsilon}}{l^2}
\int_{(\kappa+\frac 1{\log X})} |\Gamma_\delta(s)|
(1+|s|)^{1+\epsilon} |ds| \ll \frac{X^{1+2\kappa+\epsilon}}{l^2}.
$$
Using this together with (4.2) (keeping in mind that $M\le \sqrt{X}$)
we conclude that the quantity in (4.12) is bounded by
$$
\ll X^{1+\kappa +\epsilon} \Big(\frac{1}{l^2} +\frac{X^{-\text{Re }\delta_1/2}
+X^{-\text{Re }\delta_2/2}}{l^{\frac 32}} + \frac{X^{-\text{Re }\tau}}{l}\Big),
$$
which when inserted in (4.11) yields the Proposition.

\head 5.  Proof of Proposition 2.3 \endhead

\noindent Observe that
$$
\Sum_M\biggl(\L{-8d}{l} A_{\delta,\tau}(d) ;\Psi\biggr)
=  \sum_{n=1}^{\infty} \frac{r_{\delta}(n)}
{\sqrt{n}} \Sum_M \biggl(\L{-8d}{ln};F_n\biggr),
\tag{5.1}
$$
where
$$
F_n(t) =F_n(\delta,\tau;t)
=\Psi(t) W_{\delta,\tau}\L{n\pi}{8Xt}.
$$
Using the Poisson summation formula, Lemma 3.6 above, we obtain
$$
\Sum_M\biggl(\L{-8d}{ln};F_n\biggr) = \frac{1}{2ln} \L{16}{ln}
 \sum\Sb \alpha\le Y\\ (\alpha,2ln)=1\endSb \frac{\mu(\alpha)}{\alpha^2}
\sum_{k=-\infty}^{\infty} (-1)^k G_{-k}(ln) {\tilde F}_n\L{kX}{2\alpha^2
ln}.
\tag{5.2}
$$

Using this in (5.1), we deduce that
$$
\Sum_M\biggl(\L{-8d}{l}A_{\delta,\tau}(d)
;\Psi\biggr)
=  {\Cal P}(l) + {\Cal R}_{0}(l),
$$
where ${\Cal P}(l)$ is the main principal contribution
(arising from the $k=0$ term in (5.2)), and ${\Cal R}_{0}(l)$
includes all the non-zero terms $k$ in
(5.2).  Thus
$$
{\Cal P}(l) = \frac{1}{2l}\sum_{n=1}^{\infty}
\frac{r_{\delta}(n)}{n^{\frac 32}} \L{16}{ln}
\sum\Sb \alpha \le Y\\ (\alpha,2ln)
=1 \endSb \frac{\mu(\alpha)}{\alpha^2}  G_0(ln){\tilde F_n}(0),
$$
and
$$
{\Cal R}_{0}(l)
= \frac{1}{2l}\sum_{n=1}^{\infty} \frac{r_{\delta}(n)}
{ n^{3\over 2}}
\L {16}{l n} \sum\Sb \alpha \le Y\\ (\alpha, 2 l n) =1 \endSb
\frac{\mu(\alpha)}{\alpha^2}
\sum\Sb k=-\infty\\ k\neq 0\endSb^{\infty}
 (-1)^k G_{-k}(ln) {\tilde F}_n \left( \frac{kX}
{2\alpha^2 l n}\right). \tag{5.3}
$$

\subhead 5.1  The principal ${\Cal P}(l)$ contribution
\endsubhead

\noindent Note that ${\tilde F}_n(0) = {\hat F}_n(0)$ and that
$G_0(ln)=\phi(ln)$ if $ln=\square$ and $G_0(ln)=0$ otherwise.
Using this together with
$$
\sum\Sb \alpha\le Y\\ (\alpha,2ln)=1\endSb \frac{\mu(\alpha)}{\alpha^2}
= \frac{1}{\zeta(2)} \prod_{p|2ln} \Big( 1-\frac 1{p^2}\Big)^{-1}
\Big( 1+ O\Big(\frac 1Y\Big)\Big),
$$
we deduce that
$$
{\Cal P}(l)
= \frac{1+O(Y^{-1})}{\zeta(2)}
\sum\Sb n=1\\ ln=\square\endSb^{\infty}
\frac{r_{\delta}(n)}{n^{\frac12}} \L{16}{ln}
\prod_{p|2ln} \L{p}{p+1} {\hat F}_n(0).
$$

Recall that $l = l_1 l_2^2$ where $l_1$ and $l_2$ are odd, and $l_1$
is square-free.  The condition that $l n =\square$ is thus
equivalent to $ n= l_1 m^2$ for some integer $m$.
Hence
$$
\align
{\Cal P}(l)
&= \frac{1+O(Y^{-1})}{\zeta(2)\sqrt{l_1}}
\sum\Sb m=1\\ m \text{ odd}\endSb^{\infty}
\frac{r_{\delta}(l_1m^2)}{m}
\prod_{p|2lm} \L{p}{p+1} {\hat  F}_{l_1m^2}(0).\\
\endalign
$$
For any $c>|\text{Re }\tau|$ we have
$$
\align
{\hat F}_{l_1m^2} (0) &= \int_0^{\infty} \Psi(t) W_{\delta,\tau}
\L{l_1m^2 \pi}{8Xt} dt \\
&=\frac{1}{2\pi i} \int_{(c)} \Gamma_\delta(s) \L{8X}{l_1m^2 \pi}^s
\Big( \int_0^{\infty} \Psi(t) t^s dt \Big) \frac{2s}{s^2 -\tau^2 }
ds\\
&=
\frac{1}{2\pi i} \int_{(c)} \Gamma_\delta(s) {\check \Psi}(s)
\L{8X}{l_1m^2 \pi}^s \frac{2s}{s^2 -\tau^2}ds.
\\
\endalign
$$
Thus for any $c>\kappa$
$$
{\Cal P}(l) =\frac 23 \frac{1+O(Y^{-1})}{\zeta(2)\sqrt{l_1}} I(l),
\tag{5.4a}
$$
where
$$
I(l)=\frac{1}{2\pi i} \int_{(c)} \Gamma_{\delta}(s)
\L{8X}{l_1\pi}^s {\check \Psi}(s) \frac{2s}{s^2-\tau^2}
\sum\Sb m=1\\ m\text{ odd}\endSb^{\infty} \frac{r_\delta(l_1m^2)}{m^{1+2s}}
\prod_{p|lm} \L{p}{p+1} ds. \tag{5.4b}
$$

\proclaim{Lemma 5.1}  Suppose $l=l_1l_2^2$ is as above.
Then for Re $s>1 + 2|\text{Re }\delta|$
$$
\sum\Sb m=1\\ m
\text{ odd}\endSb^{\infty} \frac{r_{\delta}(l_1m^2)}{m^s}
\prod_{p|lm}
\L {p}{p+1}
= r_{\delta}(l_1) Z(s;\delta)
\eta_{\delta}(s;l)
$$
where $Z$ and $\eta$ are as defined in (2.12) and (2.13).
\endproclaim

\demo{Proof} This follows by comparing the Euler factors on both sides.
\enddemo

Using Lemma 5.1 in (5.4b), we deduce that
$$
\align
I(l)&:= \frac{r_\delta(l_1)}{2\pi i}
\int_{(c)}
\Gamma_{\delta}(s)
\L{8X}{l_1 \pi}^s {\check \Psi}(s) \frac{2s}{s^2-\tau^2}
Z(1+2s;\delta)\eta_{\delta}(1+2s;l)
ds.\\
\endalign
$$
Note first  that taking $c=\kappa+\epsilon$ here we deduce easily that
$I(l) \ll |r_\delta(l)| (X/l_1)^{\kappa+\epsilon}$.
We now
move the line of integration above to the line
Re $s=-\frac 14 +\epsilon$.  We encounter simple poles at $s=\pm \tau$, $\pm
\delta$.  The remaining integral on the $-\frac 14 +\epsilon$ line we
bound as follows:  From [3] we know that on this line
$|Z(1+2s;\delta)| \ll (1+|s|)^3$, and plainly $|\eta_{\delta}(1+2s;l)|
\ll \prod_{p|l_1}(1+O(\frac{1}{\sqrt{p}})) \prod_{p\nmid l_1}
(1+O(\frac 1{p^{1+\epsilon}})) \ll l_1^{\epsilon}$.
Hence the integral on the Re $s=-\frac 14 +\epsilon$ line is
$$
\ll \frac{|r_\delta(l_1)|l_1^{\frac 14+\epsilon}}{X^{\frac 14 -\epsilon}}
\int_{(-\frac 14 +\epsilon)}
|s|^2
|{\check \Psi}(s)| \Gamma_{\delta}(s)| |ds|
\ll  \frac{|r_\delta(l_1)|
l_1^{\frac 14+\epsilon}}{X^{\frac 14-\epsilon}}.
$$
We deduce that
$$
\align
I(l)&=  r_\delta(l_1)
\mathop{\text{Res}}_{s=\pm \delta, \pm \tau}
\biggl\{ \Gamma_{\delta}(s) \L{8X}{l_1 \pi}^s {\check \Psi}(s)
\frac{2s}{s^2-\tau^2}Z(1+2s;\delta)\eta_{\delta}(1+2s;l)
\biggr\}\\
&\hskip 2 in +
O\biggl( \frac{|r_\delta(l_1)|l_1^{\frac 14+\epsilon}}{X^{\frac 14
-\epsilon}} \biggr).\\
\endalign
$$

Using this in (5.4a), we conclude that
$$
\align
{\Cal P}(l)&=  \frac{2r_\delta(l_1) }{3\zeta(2)\sqrt{l_1}}
\mathop{\text{Res}}\Sb s=\pm \delta\\ s= \pm \tau\endSb
\biggl\{ \Gamma_{\delta}(s)\L{8X}{l_1 \pi}^s {\check \Phi}(s)
\frac{2s}{s^2-\tau^2}Z(1+2s;\delta)\eta_{\delta}(1+2s;l)
\biggr\}\\
&\hskip 2 in +
O\biggl(\frac{|r_\delta(l_1)|X^{\kappa+\epsilon}}{Y l_1^{\frac 12+\kappa}}
+ \frac{|r_\delta(l_1)|X^{\epsilon}}{(Xl_1)^{\frac 14}}
\biggr). \tag{5.5}\\
\endalign
$$

\subhead 5.2 Extracting the secondary principal term from ${\Cal R}_0(l)$
\endsubhead

\noindent Define for all real numbers $\xi$, and all complex numbers
$w$ with Re $w >0$,
$$
f(\xi,w) = \int_0^{\infty} {\tilde F}_t\L{\xi}{t} t^{w-1} dt. \tag{5.6}
$$
Since $|{\tilde F}_t(\xi/t)| \le 2|{\tilde F}_t(0)| \le e^{-\frac{t}{20X}}$
by Lemma 3.1, clearly the integral above is absolutely
convergent for Re $ w>0$.  We collect below some properties of $f(\xi,w)$
which are easily established by making minor modifications
to the proof of Lemma 5.2 of [11].

\proclaim{Lemma 5.2} For corresponding choices of sign define
$$
{\Bbb G}_{\pm}(u)= (2\pi)^{-u} \Gamma(u) \Big( \cos(\tfrac {\pi}{2} u)
\pm \sin(\tfrac{\pi}{2}u)\Big).
$$
If $\xi \neq 0$ then for any $1+\text{Re } w > c > \max(|\text{Re }\tau|,
\text{Re } w)$ we have
%$$
%\align
%f(\xi,w) &= |\xi|^w \int_0^{\infty}
%\biggl\{ {\check \Phi}(w+\tau) \L{8X}{\pi}^\tau \omega_{\delta_1,\delta_2}
%\L{\pi |\xi|}{8Xz} + {\check \Phi}(w-\tau) \L{8X}{\pi}^{-\tau}
%\omega_{-\delta_1,-\delta_2}\L{\pi|\xi|}{8Xz}\biggr\}\\
%&\hskip 1 in \times  (\cos (2\pi z)
%+ \text{sgn}(\xi) \sin(2\pi z)) \frac{dz}{z^{w+1}}.
%\tag{4.9a}\\
%\endalign
%$$
$$
f(\xi,w)= |\xi|^{w} {\check \Phi}(w)
\frac{1}{2\pi i} \int_{(c)} \Gamma_{\delta}(s) \L{8X}{\pi |\xi|}^s
{\Bbb G}_{\text{sgn}(\xi)}(s-w)\frac{2s}{s^2-\tau^2}ds.
\tag{5.7}
$$
For $\xi\neq 0$, $f(\xi,w)$ is a holomorphic function of $w$ in
$\text{Re } w > -1+|\text{Re }\tau|$, and in the region
$1\ge \text{Re }w > -1+|\text{Re }\tau|$ satisfies the bound
$$
|f(\xi,w)| \ll (1+|w|)^{-\text{Re }w -\frac 12} \exp\Big(-\frac{1}{10}
\frac{\sqrt{|\xi|}}{\sqrt{X(1+|w|)}}\Big) |\xi|^w |{\check \Phi}(w)| .
$$
\endproclaim

%Observe that for any sequence of numbers $a_n \ll n^\epsilon$, and
%any smooth function $g$ with $g(0)=0$ and $g(x)$ decaying rapidly as
%$x\to \infty$,  we have
Using the Mellin transform identity (4.9),
%$$
%\sum_{n=1}^{\infty} a_n g(n) = \frac{1}{2\pi i} \int_{(c)} \ \
%\sum_{n=1}^{\infty} \frac{a_n}{n^w} \biggl(\int_0^\infty g(t) t^{w-1}
%dt \biggr) dw,
%$$
%where $c >1$.  Hence
we may recast the expression for
${\Cal R}_0(l)$ (see (5.3) above) as
$$
{\Cal R}_0(l)=
\frac{1}{2l} \sum\Sb \alpha \le Y\\ (\alpha, 2l)=1\endSb
\frac{\mu(\alpha)}{\alpha^2}
\sum\Sb k=-\infty \\ k\neq 0\endSb^{\infty}
\frac{(-1)^k}{2\pi i}
\int_{(c)} \sum\Sb n=1 \\ (n,2\alpha)=1\endSb^{\infty}
\frac{r_{\delta}(n)}{n^{\frac 32+w} } G_{-4k}(ln)
f\biggl( \frac{kX}{2\alpha^2 l}, w\biggr) dw, \tag{5.8}
$$
for any $c> |\text{Re }\delta|$.

\proclaim{Lemma 5.3} Write $-4k=k_1k_2^2$ where $k_1$ is a fundamental
discriminant (possibly $k_1=1$, giving the trivial character),
and $k_2$ is positive.  In the region Re $s >1+ |\text{Re } \delta|$
$$
\align
\sum\Sb n=1\\ (n,2a)=1\endSb^{\infty}
\frac{r_{\delta}(n)}{n^s} \frac{G_{-4k}(ln)}{\sqrt{n}}
&= L(s-\delta,\chi_{k_1})L(s+\delta,\chi_{k_1})
\prod_{p} {\Cal G}_{\delta;p}(s;-k,l,\alpha) \\
&=: L(s-\delta,\chi_{k_1})L(s+\delta,\chi_{k_1})
{\Cal G}_{\delta}(s;-k,l,\alpha), \tag{5.9} \\
\endalign
$$
where ${\Cal G}_{\delta;p}(s;-k,l,\alpha)$ is defined as follows:  If
$p|2\alpha$ then
$$
{\Cal G}_{\delta;p}(s;-k,l,\alpha)=
\biggl( 1-\frac{1}{p^{s-\delta}}\L{k_1}{p}
\biggr)\biggl( 1-\frac{1}{p^{s+\delta}}\L{k_1}{p}
\biggr),
$$
while, if $p\nmid 2\alpha$,
$$
{\Cal G}_{\delta;p}(s;-k,l,\alpha)
= \biggl(1-\frac{1}{p^{s-\delta}}
\L{k_1}{p}\biggr)\biggl(1-\frac{1}{p^{s+\delta}} \L{k_1}{p}\biggr)
\ \
\sum_{r=0}^{\infty} \frac{r_{\delta}(p^{r})}{p^{rs}}
\frac{G_{-4k}(p^{r+\text{ord}_p(l)})}
{p^{\frac r2}}.
$$
Then ${\Cal G}_{\delta}(s;-k,l,\alpha)$
is holomorphic in the region Re $s >\frac 12+|\text{Re }\delta|$, and
for Re $s \ge \frac{1}{2} + |\text{Re } \delta| + \epsilon$
satisfies the bound
$$
| {\Cal G}_{\delta}(s;-k,l,\alpha)| \ll \alpha^\epsilon |k|^{\epsilon}
l^{\frac 12+\epsilon} (l,k_2^2)^{\frac 12}.
\tag{5.10}
$$
\endproclaim
\demo{Proof}  This follows by making minor changes to the
proof of Lemma 5.3 of [11].
\enddemo

We use Lemma 5.3 in (5.8), and move the line of integration to the
line Re $w = -\frac 12 + |\text{Re }\delta| +\epsilon$.  We
encounter poles only when $-k =\square$ (so that $k_1 =1$, and
$L(s,\chi_{k_1}) = \zeta(s)$): in this case, we have simple poles at
$w= \pm \delta$, and the residues of these poles give rise to
a second main term (see (5.12) below).  Thus we may write
${\Cal R}_0(l) = {\Cal R}(l) + {\Cal Q}(l)$ where
$$
\align
{\Cal R}(l) &= \frac{1}{2l} \sum\Sb \alpha \le Y \\ (\alpha,2l)=1\endSb
\frac{\mu(\alpha)}{\alpha^2} \sum\Sb k=-\infty \\ k\neq 0\endSb^{\infty}
\frac{(-1)^k}{2\pi i} \int_{(-\frac 12 + |\text{Re }\delta| +\epsilon)}
L(1+w+\delta,\chi_{k_1}) L(1+w-\delta,\chi_{k_1}) \\
&\hskip 2 in \times
{\Cal G}_\delta(1+w;-k,l,\alpha) f\Big( \frac{kX}{2\alpha^2 l}, w\Big) dw,
\tag{5.11}\\
\endalign
$$
and (with an obvious change of notation, writing $k^2$ in place of $-k$)
$$
{\Cal P}_2(l) = \frac{1}{2l}
\sum\Sb \alpha\le Y\\(\alpha,2l)=1\endSb \frac{\mu(\alpha)}{\alpha^2}
\sum_{\mu=\pm}  \zeta(1+2\mu \delta)
\sum_{k=1}^{\infty} (-1)^k {\Cal G}_{\delta}(1+\mu\delta;k^2,l,\alpha)
f\Big(-\frac{k^2 X}{2\alpha^2 l},\mu\delta\Big).\tag{5.12}
$$

\subhead 5.3 The secondary principal term ${\Cal Q}(l)$
\endsubhead

\noindent Define, for $\mu=\pm$,
and for any $u$ with Re $u >\frac 12$
$$
{\Cal H}_{\mu,\delta}(u;l,\alpha) = l^u \sum_{k=1}^{\infty}
\frac{(-1)^k}{k^{2u}} {\Cal G}_{\delta}(1+\mu\delta;k^2,l,\alpha).
$$
Note that the above series converges absolutely when Re $ u>\frac 12$.
Using Lemma 5.2, we see that
$$
\sum_{\mu =\pm}
\sum_{k=1}^{\infty} (-1)^k {\Cal G}_{\delta}(1+\mu\delta;k^2,l,\alpha)
f\Big(-\frac{k^2X}{2\alpha^2 l},\mu\delta\Big)
$$
may be recast as
$$
\frac{1}{2\pi i} \int_{(c)} \Gamma_{\delta}(s)
\L{16\alpha^2}{\pi}^s \sum_{\mu = \pm}
{\check \Psi}(\mu\delta) \L{X}{2\alpha^2}^{\mu\delta}
{\Bbb G}_{-}(s-\mu\delta)\frac{2s}{s^2-\tau^2}
{\Cal H}_{\mu,\delta}(s-\mu\delta;l,\alpha) ds,
\tag{5.13}
$$
where $c> \max(\frac 12 +|\text{Re }\delta|, |\text{Re }\tau|)= \frac 12+ |
\text{Re }\delta|$.

 From the definition of ${\Cal G}_\delta$ we see that
$$
\align
{\Cal H}_{\mu,\delta}(u;l,\alpha) &=
-l^{u} (1-2^{1-2u}) \sum_{k=1}^{\infty} \frac{1}{k^{2u}} {\Cal G}_{\delta}
(1+\mu\delta;k^2,l,\alpha) \\
&=
-l^u (1-2^{1-2u}) \prod_p \sum_{b=0}^{\infty} \frac{{\Cal G}_{\delta;p}
(1+\mu\delta;p^{2b},l,\alpha)}{p^{2bu}}.\\
\endalign
$$
Using the expression for ${\Cal G}_{\delta;p}$ in Lemma 5.3 and then
employing Lemma 3.3 to evaluate it, we see that we may write
$$
{\Cal H}_{\mu,\delta}(u;l,\alpha) = - l (1-2^{1-2u}) l_1^{u-\frac 12}
\zeta(2u) \zeta(2u+1+4\mu\delta) {\Cal H}_{1,\mu,\delta}(u;l,\alpha)
$$
where ${\Cal H}_1 =\prod_p {\Cal H}_{1;p}$ with
$$
{\Cal H}_{1;p} =
\cases
\Big(1-\frac 1p\Big) \Big(1-\frac{1}{p^{1+2\mu\delta}}\Big)\Big(1-\frac{1}
{p^{2u+1+4\mu\delta}}\Big) &\text{if } p|2\alpha \\
\Big(1-\frac{1}{p}\Big) \Big(1-\frac{1}{p^{1+2\mu\delta}}\Big)
\Big(1+\frac 1p +\frac 1{p^{1+2\mu\delta}} -\frac{1}{p^{2u+2+4\mu\delta}}
\Big) &\text{if } p\nmid 2\alpha l\\
\Big(1-\frac{1}{p}\Big) \Big(1-\frac{1}{p^{1+2\mu\delta}}\Big)
\Big(1+\frac{1}{p^{2u+2\mu\delta}}\Big) &\text{if } p|l_1\\
\Big(1-\frac{1}{p}\Big) \Big(1-\frac{1}{p^{1+2\mu\delta}}\Big)
\Big(1+\frac{1}{p^{1+2\mu\delta}}\Big) &\text{if } p|l, p\nmid l_1.\\
\endcases
$$
These expressions show that ${\Cal H}_{\mu,\delta}(u;l,\alpha)$
is analytic in the domain Re $ u > -\frac 12 + |\text{Re }\delta|$.
Thus we may move the line of integration in (5.13) to the
line Re $s = \kappa+ \frac{1}{\log X}$, and
since we encounter no poles, (5.13) is given by the
resulting integral on this line.  Using these observations in
(4.12) we conclude that
$$
\align
{\Cal P}_2(l)  &=\frac{1}{2\pi i} \int_{(\kappa+\frac{1}{\log X})}
\Gamma_\delta(s) \L{16}{\pi}^s \frac{2s}{s^2-\tau^2}
 \sum_{\mu=\pm} {\check \Psi}(\mu\delta)
\L{X}{2}^{\mu\delta} \zeta(1+2\mu\delta)
 {\Bbb G}_-(s-\mu\delta) \\
&\hskip 1.3 in \times
\frac{1}{2l}\sum\Sb \alpha \le Y\\(\alpha,2l)=1\endSb \frac{\mu(\alpha)}
{\alpha^{2-2s+2\mu\delta}}
{\Cal H}_{\mu,\delta}(s-\mu\delta;l,\alpha) ds.\tag{5.14}\\
\endalign
$$

We now wish to show that the sum over $\alpha$ in (5.14) may
be extended to infinity, at the cost of an acceptable error.
Let ${\Cal C}$ denote the closed curve oriented counter-clockwise
consisting of the following four line segments: from
$|\text{Re }\delta| +1/(2\log  X) - i(|\text{Im } \delta|+
 1/\log  X)$
to $|\text{Re }\delta| +1/(2\log X) + i(|\text{Im } \delta|+1/\log X)$,
and from there to
$-|\text{Re }\delta| -1/(2\log X) + i(|\text{Im } \delta|+1/\log X)$, and
from here to
$-|\text{Re }\delta| -1/(2\log X) - i(|\text{Im } \delta|+1/\log X)$
and lastly back to
$|\text{Re }\delta| +1/(2\log X) - i(|\text{Im } \delta|+1/\log X)$.
Given $s$ with Re $s =\kappa +\frac{1}{\log X}$, the
 function $2z {\check \Psi}(z) \L{X}{2\alpha^2}^z \zeta(1+2z) {\Bbb G}_-
(s-z) {\Cal H}_{\mu,\delta}(s-z;l,\alpha)$ is analytic for $z$ inside
 ${\Cal C}$.  So by Cauchy's theorem
$$
\align
&\sum_{\mu=\pm} {\check \Psi}(\mu\delta)
\L{X}{2\alpha^2}^{\mu\delta} \zeta(1+2\mu\delta)
 {\Bbb G}_-(s-\mu\delta)
{\Cal H}_{\mu,\delta}(s-\mu\delta;l,\alpha)\\
=& \frac{1}{2\pi i} \int_{{\Cal C}}
{\check \Psi}(z) \L{X}{2\alpha^2}^z \zeta(1+2z) {\Bbb G}_-
(s-z) {\Cal H}_{\mu,\delta}(s-z;l,\alpha) \frac{2z}{z^2 -\delta^2} dz.
\tag{5.15}\\
\endalign
$$
We now assume that $|\text{Im }z| \le 1$, say.  For $z$ on ${\Cal C}$
we see that $\kappa + |\text{Re }\delta| + 3/(2\log X)
 \ge \text{Re }(s-z) \ge 1/(2\log X)$.
Further $|\frac{2z}{z^2-\delta^2}{\check \Psi}(z) \zeta(1+2z)
\L{X}{2\alpha^2}^z|
\ll (\log^2 X) (X\alpha^2)^{|\text{Re }\delta|+1/\log X}$, and by Stirling's
formula we see
that $|{\Bbb G}_-(s-z)| \ll (1+|\text{Im }(s)|)^{2\kappa -1/2} \ll 1$.
Lastly from our expressions for ${\Cal H}_{\mu,\delta}(u;l,\alpha)$
we deduce that $|{\Cal H}_{\mu,\delta} (u;l,\alpha) |
\ll l^{1+\epsilon} l_1^{\text{Re }u-\frac 12} \alpha^{\epsilon} (1+|s|)$.
From these estimates we conclude that (5.15) is
bounded by
$$
l^{1+\epsilon} l_1^{\kappa +|\text{Re }\delta|-\frac 12}
 (X\alpha^2)^{|\text{Re }\delta|+\epsilon} (1+|s|).
$$
We deduce that
$$
\frac{1}{2l} \sum\Sb \alpha>Y \\ (\alpha,2l)=1\endSb \frac{\mu(\alpha)}
{\alpha^{2-2s}}
\sum_{\mu=\pm} {\check \Psi}(\mu\delta)
\L{X}{2\alpha^2}^{\mu\delta} \zeta(1+2\mu\delta)
 {\Bbb G}_-(s-\mu\delta)
{\Cal H}_{\mu,\delta}(s-\mu\delta;l,\alpha)
$$
is bounded by
$\ll l^{\epsilon} l_1^{\kappa+|\text{Re } \delta|
-\frac 12} X^{|\text{Re }\delta|+\epsilon}
(1+|s|)
Y^{-1 +2 \kappa +2|\text{Re }\delta|}$.
Using this in (5.14) we conclude that the error incurred in
extending the sum over $\alpha$ to infinity is
$$
\align
&\ll l^{\epsilon} l_1^{\kappa +|\text{Re }\delta|
-\frac 12} X^{|\text{Re }\delta|+\epsilon}
Y^{-1+2\kappa +2|\text{Re }\delta|}
\int_{(\kappa)} |\Gamma_\delta(s)| (1+|s|)
\frac{|s|}{|s^2-\tau^2|}|ds|
\\
&\ll l^{\epsilon} l_1^{\kappa+|\text{Re } \delta|
-\frac 12} X^{|\text{Re }\delta|+\epsilon}
Y^{-1+2\kappa +2|\text{Re }\delta|}.
\\
\endalign
$$

%From our expression for ${\Cal H}_{\mu,\delta}(u;l,\alpha)$ we
%deduce that when Re$(u) \ge \frac{1}{\log X}$ say
%$$
%|{\Cal H}_{\mu,\delta}(u;l,\alpha)| \ll l^{1+\epsilon}
%l_1^{\text{Re}(u)-\frac 12} (|u|+1) X^{\epsilon}.
%$$

%We now extend the sum over $\alpha$ to $\infty$; the
%error incurred in doing so is negligible.
Thus, up to an error
$O(l^{\epsilon} l_1^{\kappa+|\text{Re }\delta|
-\frac 12} X^{|\text{Re }\delta|+
\epsilon}Y^{-1+2\kappa +2|\text{Re }\delta|})$,
${\Cal P}_2(l)$ is given by
$$
\sum_{\mu=\pm} {\check \Psi}(\mu\delta)
\L{X}{2}^{\mu\delta} \zeta(1+2\mu\delta)
\frac{1}{2\pi i} \int_{(\kappa)}
\Gamma_\delta(s) \L{16}{\pi}^s {\Bbb G}_-(s-\mu\delta) \frac{2s}{s^2-\tau^2}
{\Cal K}_{\mu,\delta}(s;l) ds,
\tag{5.16}
$$
where
$$
{\Cal K}_{\mu,\delta}(s;l)= \frac{1}{2l}
\sum\Sb \alpha =1\\(\alpha,2l)=1\endSb^{\infty}
 \frac{\mu(\alpha)}{\alpha^{2-2s+2\mu\delta}}
{\Cal H}_{\mu,\delta}(s-\mu\delta;l,\alpha).
$$

Using our expression for ${\Cal H}_{\mu,\delta}$ a calculation
gives
$$
\align
{\Cal K}_{\mu,\delta}(s;l) &=
-\frac 1{4l_1^{\frac 12+\mu\delta}} \frac {\phi(l)}{l}
\prod_{p|2l} \biggl( 1-\frac{1}{p^{1+2\mu\delta}}\biggr)
\prod\Sb p|l\\p\nmid l_1\endSb \biggl( 1+ \frac{1}{p^{1+2\mu \delta}}\biggr)
r_{s}(l_1) \\
&\times
\biggl(\frac{4^s+4^{-s}-2^{-1-2\mu\delta}-2^{1+2\mu\delta}}{4^s}\biggr)
\zeta(2s-2\mu\delta) \zeta(2s+1+2\mu\delta) \\
&\times
\prod_{p\nmid 2l} \biggl (1-\frac 1p\biggr)\biggl( 1-
\frac{1}{p^{1+2\mu\delta}}\biggr)\biggl(1+\frac 1p +\frac1{p^{1+2\mu\delta}}
+ \frac{1}{p^{3+4\mu\delta}}-\frac{(p^{2s}+p^{-2s})}{p^{2+2\mu\delta}}
\biggr).\\
\endalign
$$
Using this together with the
functional equation for $\zeta(s)$ and the relations
$\Gamma(z) \Gamma(1-z) = \pi \text{cosec}(\pi z)$ and
$\Gamma(z) \Gamma(z+\frac 12)
= \pi^{\frac12} 2^{1-2z} \Gamma(2z)$ we see that
$$
{\Cal J}_{\mu,\delta}(s;l) := \Gamma_\delta(s) {\Bbb G}_-(s-\mu\delta)
\L{16}{\pi}^s {\Cal K}_{\mu,\delta}(s;l)
$$
satisfies the functional equation ${\Cal J}_{\mu,\delta}(s;l) =
{\Cal J}_{\mu,\delta}(-s;l)$.  In fact, we obtain the useful identity
$$
\zeta(1+2\mu \delta) {\Cal J}_{\mu,\delta} (s;l)
= \frac{2r_{s}(l_1)}{3\zeta(2) \sqrt{l_1}}\L{16}{\pi l_1}^{\mu\delta}
\Gamma_s(\mu\delta) Z(1+2\mu\delta;s) \eta_s(1+2\mu\delta;l);
\tag{5.17}
$$
it is plain that the left side above is invariant under $s\to -s$.

Consider now for $\mu =\pm$ the integral in (5.16): that is,
$$
\frac{1}{2\pi i} \int_{(\kappa)}
{\Cal J}_{\mu,\delta}(s;l) \frac{2s}{s^2-\tau^2} ds.
\tag{5.18}
$$
We move the line of integration to the line Re$(s) = - \kappa
$.   We encounter simple poles at $s=\delta, -\delta,
\tau$, and $-\tau$.  Thus (5.18) equals
$$
\mathop{\text{Res}}_{s=\pm\delta, \pm\tau} {\Cal J}_{\mu,\delta} (s;l)
\frac{2s}{s^2-\tau^2}
+ \frac{1}{2\pi i} \int_{(-\kappa)}
{\Cal J}_{\mu,\delta}(s;l) \frac{2s}{s^2-\tau^2} ds.
$$
Changing $s$ to $-s$ and using the relation ${\Cal J}_{\mu,\delta}(s;l)
= {\Cal J}_{\mu;\delta}(-s;l)$ we
see that the above is
$$
= \mathop {\text{Res}}_{s=\pm\delta, \pm\tau} {\Cal J}_{\mu,\delta} (s;l)
\frac{2s}{s^2-\tau^2}
- \frac{1}{2\pi i} \int_{(\kappa)}
{\Cal J}_{\mu,\delta}(s;l) \frac{2s}{s^2-\tau^2} ds.
$$
Hence (5.18) equals
$$
\frac{1}{2} \biggl( \mathop{\text{Res}}_{s=\pm\delta}
{\Cal J}_{\mu,\delta} (s;l)
\frac{2s}{s^2-\tau^2} + {\Cal J}_{\mu,\delta}(\tau;l) + {\Cal J}_{\mu,\delta}
(-\tau;l) \biggr)
= \mathop{\text{Res}}_{s=\mu\delta} {\Cal J}_{\mu,\delta} (s;l)
\frac{2s}{s^2-\tau^2} + {\Cal J}_{\mu,\delta}(\tau;l),
$$
using once again that ${\Cal J}_{\mu,\delta}(s;l)
={\Cal J}_{\mu,\delta}(-s;l)$.

We conclude that
$$
\align
{\Cal P}_2(l) &= \sum_{\mu=+,-} {\check \Psi}(\mu\delta)
\L{X}{2}^{\mu\delta} \zeta(1+2\mu\delta) \left(
 \mathop{\text{Res}}_{s=\mu\delta}
{\Cal J}_{\mu,\delta} (s;l)
\frac{2s}{s^2-\tau^2} + {\Cal J}_{\mu,\delta}(\tau;l)
\right) \\
&\hskip 1.5 in
+ O\Big(\frac{l^\epsilon X^{|\text{Re }\delta| +\epsilon}
 l_1^{\kappa+|\text{Re }\delta| -\frac 12}}{Y^{1-2\kappa -2|\text{Re }
\delta|}}
\Big). \tag{5.19}
\\
\endalign
$$

\subhead 5.4 The contribution of the remainder terms
${\Cal R}(l)$ \endsubhead

\noindent The contribution of the remainder terms ${\Cal R}(l)$ is
bounded in much the same manner as the analogous quantity in [11] (see
Section 5.4 there).  For the sake of completeness we give a detailed sketch
of the main ideas of the proof.

First we bound $|{\Cal R}(l)|$
for individual $l$.  Using the bounds of Lemmas
5.2 and 5.3 in (5.11) we get that $|{\Cal R}(l)|$ is
$$
\align
&\ll \frac{l^{-|\text{Re }\delta|+\epsilon}}{X^{\frac 12-|\text{Re }\delta|
-\epsilon}} \sum_{\alpha \le Y} \frac{1}{\alpha^{1+2|\text{Re }\delta|
-\epsilon}} \int_{(-\frac 12+ |\text{Re }\delta| +\epsilon)}
|{\check \Psi} (w)| (1+|w|)^{-|\text{Re }\delta|}
\\
&\hskip .25 in \times \sum\Sb k=-\infty\\k\neq 0\endSb^{\infty}
\frac{|L(1+w+\delta,\chi_{k_1})L(1+w-\delta,\chi_{k_1})|}
{|k_1|^{\frac 12-|\text{Re }\delta|}}
k_2^{2|\text{Re }\delta|} \exp\Big( -\frac 1{10} \frac{\sqrt{|k|}}
{\alpha\sqrt{l(1+|w|)}}\Big) |dw|.
\\
\endalign
$$
Performing the sum over $k_2$ we see that this is bounded by
$$
\align
&\frac{l^{\frac 12+\epsilon}}{X^{\frac 12-|\text{Re }\delta|
-\epsilon}} \sum_{\alpha \le Y} \alpha^{\epsilon}
\int_{(-\frac 12+ |\text{Re }\delta| +\epsilon)}
|{\check \Psi} (w)| (1+|w|)^{\frac 12}
\\
&\hskip .9 in\times  \sum_{k_1}
\frac{|L(1+w+\delta,\chi_{k_1})L(1+w-\delta,\chi_{k_1})|}
{|k_1|}
\exp\Big( -\frac 1{10} \frac{\sqrt{|k_1|}}
{\alpha\sqrt{l(1+|w|)}}\Big) |dw|.
\\
\endalign
$$
We split the $k_1$ into dyadic blocks and use Cauchy's inequality
with Lemma 3.5 to estimate
these contributions.  We deduce that
$$
\align
|{\Cal R}(l)| &\ll
\frac{l^{\frac 12+\epsilon} Y^{1+\epsilon}}{X^{\frac 12 -|\text{Re }\delta|
-\epsilon}} \int_{(-\frac 12+ |\text{Re }\delta| +\epsilon)}
|{\check \Psi} (w)| (1+|w|) |dw|
\ll
 \frac{l^{\frac 12+\epsilon} Y^{1+\epsilon}}{X^{\frac 12 -|\text{Re }\delta|
-\epsilon}} \Psi_{(2)}\Psi_{(3)}^{\epsilon}.\\
\endalign
$$

We now sketch how a better bound for ${\Cal R}(l)$ may be obtained
on average.
Let $\beta_l = \frac{\overline{{\Cal R}(l)}}{|{\Cal R}(l)|}$ if ${\Cal R}(l)
\neq 0$, and $\beta_l =1$ otherwise.  Then, from (5.11),
$\sum_{l=L}^{2L-1} |{\Cal R}(l)| =\sum_{l=L}^{2L-1} \beta_l {\Cal R}(l)$
is
$$
\align
\ll& \sum\Sb \alpha \le Y\\ (\alpha,2)=1\endSb \frac{1}{\alpha^2}
 \int_{(-\frac 12+|\text{Re }
\delta|+\epsilon)} \sum\Sb k=-\infty \\ k\neq 0\endSb^{\infty}
|L(1+w+\delta,\chi_{k_1})L(1+w-\delta,\chi_{k_1})| \\
&\hskip .9 in \times \biggl|\sum\Sb l=L\\(l,\alpha)=1\endSb^{2L-1}
\frac{\beta_l}{l}
{\Cal G}_\delta(1+w;-k,l,\alpha) f\biggl(\frac{kX}{2\alpha^2l},w\biggr)
\biggr| |dw|. \tag{5.20}\\
\endalign
$$
We now split the sum over $k$ into dyadic blocks $K\le |k|\le 2K-1$.
By Cauchy's inequality the sum over $k$ in this range is
bounded by the product of two terms.  The first of these terms
is
$$
\Big(\sum_{|k|=K}^{2K-1} k_2
|L(1+w+\delta,\chi_{k_1})L(1+w-\delta,\chi_{k_1})|^2
\Big)^{\frac 12}
\ll (K(1+|w|))^{\frac 12+\epsilon},
$$
upon using Cauchy's inequality again with Lemma 3.5.
The second term in question is
$$
\biggl(\sum_{|k|=K}^{2K-1} \frac{1}{k_2}
\biggl|\sum\Sb l=L\\
(l,2\alpha)=1\endSb^{2L-1}
\frac{\beta_l}{l} {\Cal G}_\delta(1+w;-k,l,\alpha)
f\biggl(\frac{kX}{2\alpha^2 l},w
\biggr)\biggr|^2\ \biggr)^{\frac 12}. \tag{5.21}
$$

\proclaim{Lemma 5.4}  Let $\alpha \le Y$, $K$ and $L$ be positive integers,
and suppose $w$ is a complex number with Re $w=-\frac 12+|\text{Re }\delta|
+\epsilon$.  Then
for any choice of complex numbers $\gamma_l$ with $|\gamma_l|\le 1$ we
have
$$
\sum_{|k|=K}^{2K-1} \frac{1}{k_2} \
\biggl| \sum\Sb l=L\\ (l,2\alpha)=1\endSb^{2L-1}
\frac{\gamma_l}{l}
{\Cal G}(1+w;-k,l,\alpha)
f\biggl(\frac{kX}{2\alpha^2 l},w\biggr) \biggr|^2
$$
is bounded by
$$
(1+|w|)^{-2|\text{Re }\delta|+\epsilon} |{\check \Psi}(w)|^2
\frac{\alpha^{2-4|\text{Re }\delta|+\epsilon} L^{2-2|\text{Re }\delta|
+\epsilon} K^{2|\text{Re }\delta|+ \epsilon}}{X^{1-2|\text{Re }\delta| -
\epsilon}}
\exp\left(-\frac{1}{20} \frac{\sqrt{K}}{\alpha \sqrt{L (1+|w|)}}\right),
$$
and also by
$$
((1+|w|)\alpha KLX)^{\epsilon} |{\check \Psi}(w)|^2
\Big(\frac{\alpha^2 L(1+|w|)}{K}\Big)^{2|\text{Re }\tau|-2|\text{Re }
\delta|}
\frac{\alpha^{2}L}
{K X^{1-2|\text{Re }\delta|}}(K+L).
$$
\endproclaim

We bound (5.21) using the first bound of the Lemma
for $K\ge \alpha^2 L(1+|w|) \log^2 X$, and the second bound for smaller
$K$.  Inserting this bound in (5.20) gives (with a little calculation)
$$
\sum_{l=L}^{2L-1} |{\Cal R}(l)| \ll
\frac{L^{1+\epsilon} }{X^{\frac 12-|\text{Re }\delta| -\epsilon}}
\sum_{\alpha \le Y}
\alpha^{\epsilon} \int_{(-\frac 12+\epsilon)} |{\check \Psi}(w)|
(1+|w|)^{1+\epsilon} |dw| \ll \frac{L^{1+\epsilon} Y^{1+\epsilon}}
{X^{\frac 12-|\text{Re }\delta| -\epsilon}}
\Psi_{(2)} \Psi_{(3)}^{\epsilon},
$$
as desired.

\demo{Proof of Lemma 5.4} We follow closely the proof of Lemma 5.4 in [11].
Using the bound for ${\Cal G}_\delta$ in Lemma 5.3, and the
bound for $|f(\xi,w)|$ in Lemma 5.2 we easily
obtain the first bound claimed.

Write the integral in (5.7) as $\frac{1}{2\pi i} \int_{(c)}
g(s,w;\text{sgn}(\xi)) \L{8X}{|\xi| \pi}^{s} ds$.  Taking $c=|\text{Re }\tau|
+ \epsilon$,
we see that (for $K\le |k|\le 2K-1$)
$$
\align
\biggl| \sum\Sb l=L\\ (l,2\alpha)=1\endSb^{2L-1}
\frac{\gamma_l}{l}{\Cal G}_\delta
(1+w;-k,l,\alpha)
&f\biggl(\frac{kX}{2\alpha^2 l} ,w\biggr) \biggr|
\ll |{\check \Psi}(w)| \frac{\alpha^{1+2|\text{Re }\tau|
-2|\text{Re }\delta|+\epsilon}}{K^{\frac 12 +|\text{Re }\tau|-
|\text{Re }\delta|-\epsilon}X^{\frac 12-|\text{Re }\delta|-\epsilon}}
\\
&\times  \int_{(c)} \biggl|g(s,w;\text{sgn}(k))
 \sum\Sb l =L\\ (l,2\alpha)=1\endSb^{2L-1} \frac{\gamma_l}{l^{1+w-s}}
{\Cal G}_\delta(1+w;-k,l,\alpha)ds \biggr|.\\
\endalign
$$
Since $|g(s,w;\text{sgn}(k))| \ll
(1+|w|)^{c-\frac 12-\text{Re }w+\epsilon}
\exp(-\frac \pi{2}
|\text{Im}(s)|)$ by Stirling's formula,
we get by Cauchy's inequality that the above is
$$
\align
&\ll (1+|w|)^{|\text{Re }\tau|-|\text{Re }\delta|+\epsilon}|{\check \Psi}(w)|
\frac{\alpha^{1+2|\text{Re }\tau|-2|\text{Re }\delta| +
\epsilon}}{K^{\frac 12-|\text{Re }\delta| +|\text{Re }\tau|
-\epsilon}X^{\frac 12-|\text{Re }\delta|-\epsilon}}
\\
&\hskip 1 in
\times \biggl( \int_{(c)} \exp\left(-\tfrac \pi{2} |\text{Im}(s)|\right)
\biggl| \sum\Sb l=L\\ (l,2\alpha)=1\endSb^{2L-1} \frac{\gamma_l}{l^{1+w-s}}
{\Cal G}_\delta (1+w;-k,l,\alpha) \biggr|^2 |ds|\biggr)^{\frac 12}.
\\
\endalign
$$
The second bound of the Lemma follows by combining this with Lemma 5.5
below.

\enddemo

\proclaim{Lemma 5.5}  Let $|\delta_l| \ll l^{\epsilon}$ be any sequence of
complex numbers, with $\delta_l=0$ if $(l,2\alpha)\neq 1$.
Let $w$ be any complex number with Re$(w)=-\frac 12
+|\text{Re } \delta| + \epsilon$.  Then
$$
\sum_{|k|=K}^{2K-1} \frac{1}{k_2}
\biggl| \sum\Sb l=L\endSb^{2L-1}
\frac{\delta_l}{\sqrt{l}} {\Cal G}_\delta (1+w;-k,l,\alpha)\biggr|^2 \ll
(KL\alpha)^{\epsilon}(K+L)L.
$$
\endproclaim

\demo{Proof}  For any integer $k =\pm \prod_{i,\ a_i\ge 1} p_i^{a_i}$
we define $a(k) = \prod_{i} p_i^{a_i+1}$, and put $b(k)=\prod_{i,\
a_i=1} p_i \ \prod_{i, \ a_i\ge 2}
p_i^{a_i-1}$.  Note that ${\Cal G}_\delta (1+w;-k,l,\alpha)=0$
unless $l$ can be written as $dm$ where $d|a(k)$ and $(m,k)=1$
with $m$ square-free.  From
the definition of ${\Cal G}$ in Lemma 5.3,  and using Lemma 3.3, we get
$$
{\Cal G}_\delta(1+w;k,l,\alpha) = \sqrt{m}\L{-k}{m} \prod_{p|m}
\biggl(1+\frac{r_\delta(p)}{p^{1+w}}\L{-k}{p}\biggr)^{-1} {\Cal G}_\delta
(1+w;-k,d,\alpha).
$$
Using Lemma 5.3 to bound $|{\Cal G}_\delta (1+w;-k,d,\alpha)|$ we see
that our desired sum is
$$
\ll (KL\alpha)^{\epsilon} \sum_{|k|=K}^{2K-1}\frac 1{k_2}
\  \sum_{d|a(k)} d
\ \biggl|\sum_{m=L/d}^{2L/d} \delta_{dm} \mu(m)^2 \L{-k}{m}
\prod_{p|m} \biggl(1+\frac{r_\delta(p)}{p^{1+w}}\L{-k}{p}\biggr)^{-1}
\biggr|^2.
$$
We interchange the sums over $d$ and $k$.  Note that $d|a(k)$ implies that
that $b(d) |k$, so that $k=b(d) f$ for some integer
$f$ with $K/b(d) \le |f| \le 2K/b(d)$.  Write $-4f=f_1f_2^2$ where
$f_1$ is a fundamental discriminant, and $f_2$ is positive.  Notice that
$k_2\ge f_2$.   Thus our desired sum is bounded
by
$$
 (KL\alpha)^{\epsilon} \sum_{d\le 2L} d \
\sum_{f =K/b(d)}^{2K/b(d)} \frac 1{f_2}
\ \biggl|\sum_{m=L/d}^{2L/d} \delta_{dm} \mu(m)^2
 \L{-fb(d)}{m} \prod_{p|m} \biggl(1+\frac{r_\delta(p)}{p^{1+w}}
\L{-fb(d)}{p}\biggr)^{-1}
\biggr|^2,
$$
and by Lemma 3.4 this is
$$
\ll (KL\alpha)^{\epsilon} \sum_{d\le L} d \frac{L}{d} \biggl(\frac{K}{b(d)}
+\frac Ld\biggr) \ll (KL\alpha)^{\epsilon}(KL+L^2) \sum_{d\le L}
\frac{1}{b(d)} \ll (KL\alpha)^{\epsilon} (KL+L^2).
$$
\enddemo

\subhead 5.5 Completion of the Proof \endsubhead

\noindent From our work above the remainder terms are under control;
and we need only simplify the main term ${\Cal P}(l) + {\Cal P}_2(l)$
arising from (5.5) and (5.19).  Using (5.17) it is
easy to see that the contribution to (5.19) from the poles at $\mu\delta$
cancel precisely the contribution to (5.5) from the poles at $\mu\delta$.
Thus our main term includes only the contribution from the
poles at $\pm \tau$ in both these expressions.  Employing (5.17) we
deduce that the main term is
$$
\align
\frac{2}{3\zeta(2)\sqrt{l_1}}
\sum_{\mu=\pm} &\biggl( r_\delta(l_1) \Gamma_\delta(\mu\tau)
\L{8X}{l_1 \pi}^{\mu\tau} {\check \Psi}(\mu\tau) Z(1+2\mu\tau;\delta)
\eta_{\delta}(1+2\mu\tau;l) \\
&+ r_{\tau}(l_1) \Gamma_\tau(\mu\delta) \L{8X}{l_1 \pi}^{\mu\delta}
{\check \Psi}(\mu\delta) Z(1+2\mu\delta;\tau) \eta_{\tau}(1+2\mu\delta;l)
\biggr).\\
\endalign
$$
This proves Proposition 2.3.

\head 6. Mollification near $s=\frac 12$:  Proof of Proposition 2.4
 \endhead

\noindent %We shall choose our mollifier to be of the form
%$$
%M(s,d) = \sum_{m} \frac{\lam(m)}{m^s} \chi_{-8d}(m)
%$$
%where $\lam(m)=0$ unless $m$ is an odd square-free integer below $M$, and
%in this case $\lam(m) = \mu(m) Q(\frac{\log M/m}{\log M})$ for
%a real-valued function $Q$ to be prescribed later.  We seek to evaluate
%${\Cal S} (|L(\tfrac 12 +\delta_1,\chi_{-8d})
%M(\tfrac 12+\delta_1,d)|^2;\Phi)$.  Set $\delta_2 =
%\overline{\delta_1}$,  so that $\tau = \text{Re }\delta_1$ is
%real, and $\delta=i\text{Im }\delta_1$ is purely imaginary.  In the following
%we shall assume that $\tau \ge - A/\log X$ for some absolute
%positive constant $A$.
From Lemma 3.2, and the definition of $\xi(s,\chi_{-8d})$ we see that
$$
\align
&{\Cal W}(\delta_1, \Phi)
= \frac{(8X/\pi)^{-\tau}}{\Gamma_\delta(\tau)}
{\Cal S}(A_{\delta,\tau}(d)|M(\tfrac 12+\delta_1,d)|^2; \Phi_{-\tau})
\\
\endalign
$$
where $\Phi_{-\tau}(t) = t^{-\tau} \Phi(t)$.
For a parameter $Y$ to be fixed later we decompose the above
as
$$
\frac{(8X/\pi)^{-\tau}}{\Gamma_\delta(\tau)}
\Big\{
{\Cal S}_M
(A_{\delta,\tau}(d)|M(\tfrac 12+\delta_1,d)|^2; \Phi_{-\tau})
+O({\Cal S}_R
(A_{\delta,\tau}(d)|M(\tfrac 12+\delta_1,d)|^2; \Phi_{-\tau}))
\Big\}.
$$
Applying Proposition 2.2 we conclude that
$$
{\Cal W}(\delta_1, \Phi)
= \frac{(8X/\pi)^{-\tau}}{\Gamma_\delta(\tau)}
{\Cal S}_M
(A_{\delta,\tau}(d)|M(\tfrac 12+\delta_1,d)|^2; \Phi_{-\tau})
+ O\Big(X^{\epsilon} \Big(\frac 1Y +X^{-\tau}\Big)\Big).
\tag{6.1}
$$

Now
$$
{\Cal S}_M(A_{\delta,\tau}(d)|M(\tfrac 12+\delta_1,d)|^2; \Phi_{-\tau})
=
\sum_{l} \Big(\sum_{rs=l}
\frac{\lam(r) \lam(s)}{r^{\frac 12+\delta_1}s^{\frac 12+\delta_2} }
\Big) {\Cal S}_M \Big(A_{\delta,\tau}(d)\L{-8d}{l};\Phi_{-\tau}\Big),
$$
and we use Proposition 2.3 to evaluate these terms.
First we note that the
various remainder terms in Proposition 2.3 contribute
(using $|\lam(n)|\ll n^\epsilon$, $r_\delta(n)\ll n^{\epsilon}$,
$d(n) \ll n^{\epsilon}$, and $\tau \ge -A/\log X$)
$$
\align
&\ll \sum_{l\le M^2} \frac{l^{\epsilon}}{l^{\frac 12+\tau}}
\Big( \frac{X^{\tau+\epsilon}}{Y l_1^{\frac 12+ \tau}} +
\frac{X^{\epsilon}}{(Xl_1)^{\frac 14}} + \frac{X^{\epsilon}
l_1^{\tau-\frac 12}}{Y^{1-2\tau}} + |{\Cal R}(l)|\Big)
\\
&\ll \frac{X^{\tau+\epsilon}}{Y} + \frac{M^{\frac 12-2\tau}X^{\epsilon}}
{X^{\frac 14}} + \frac{X^{\epsilon}}{Y^{1-2\tau}}
+ \frac{X^{\epsilon} M^{1-2\tau}Y}{\sqrt{X}} \Phi_{(2)} \Phi_{(3)}^{
\epsilon}.\\
\endalign
$$
We choose $Y=X^{\tau+\epsilon}$, and
recall that $M=X^{\frac 12-\epsilon}$.
Thus from the above considerations we get that
$$
{\Cal W}(\delta_1, \Phi)
= \frac{(8X/\pi)^{-\tau}}{\Gamma_\delta(\tau)}
\sum_{l} \Big(\sum_{rs=l}
\frac{\lam(r) \lam(s)}{r^{\frac 12+\delta_1}s^{\frac 12+\delta_2} }
\Big) {\Cal M}(l)
 + O(X^{-\tau -\epsilon} \Phi_{(2)} \Phi_{(3)}^{\epsilon}),
\tag{6.2a}
$$
where ${\Cal M}(l) = {\Cal M}_1(l)+{\Cal M}_2(l)$ with
$$
{\Cal M}_1(l) = \frac{2}{3\zeta(2)\sqrt{l_1}}
\sum_{\mu=\pm} r_\delta(l_1) \Gamma_\delta(\mu\tau)
\L{8X}{l_1 \pi}^{\mu\tau} {\check \Phi}(\mu\tau-\tau) Z(1+2\mu\tau;\delta)
\eta_{\delta}(1+2\mu\tau;l), \tag{6.2b}
$$
and
$$
{\Cal M}_2(l)
= \frac{2}{3\zeta(2)\sqrt{l_1}}
\sum_{\mu=\pm}
r_{\tau}(l_1) \Gamma_\tau(\mu\delta) \L{8X}{l_1 \pi}^{\mu\delta}
{\check \Phi}(\mu\delta-\tau)
Z(1+2\mu\delta;\tau) \eta_{\tau}(1+2\mu\delta;l). \tag{6.2c}
$$

Recall that we made the simplifying assumptions that
$\frac 14 \ge \text{Re }\delta_1 \ge -\frac{1}{\epsilon \log X} $
and that $|\delta_1| \ge \frac{\epsilon}{\log X}$.
Thus $\tau^2 - \delta^2 = |\delta_1|^2
\ge \epsilon^2/(\log X)^2$.  These assumptions
 enable us to evaluate the ${\Cal M}_1(l) $ and
${\Cal M}_2(l)$ contributions to  (6.2a) separately.
Let ${\Cal C}$ denote a closed contour (oriented counter-clockwise)
which contains the points $\pm \tau$ and such that for $w \in {\Cal C}$
we have $|\text{Re } w| \le |\tau| +C/\log X$, and $|\text{Im }w| \le
C/\log X$ for
some absolute constant $C$, and such that
$|w^2-\tau^2| \ge \epsilon^2/(3\log^2 X)$, and $|w^2-\delta^2| \ge
\epsilon^2/(3\log^2 X)$, and finally such that the perimeter
length of ${\Cal C}$ is $\ll |\delta_1|$.
Then the contribution of
${\Cal M}_1(l)$ to (6.2a) is
$$
\align
&\frac{1}{2\pi i} \int_{{\Cal C}}
\frac{2{\check\Phi}(w-\tau)}{3\zeta(2)\Gamma_\delta(\tau)}
\L{8X}{\pi}^{w-\tau}
Z(1+2w;\delta) \Gamma_\delta(w) \frac{2w}{w^2-\tau^2}
\\
&\hskip 1 in \times \biggl\{ \sum_{l}
\frac{r_{\delta}(l_1)}{l_1^{\frac 12 +w}}
\Big( \sum_{rs=l} \frac{\lam(r)\lam(s)}{r^{\frac 12+\delta_1}s^{\frac 12
+\delta_2}} \Big) \eta_{\delta}(1+2w\tau;l)\biggr\} dw.
\tag{6.3}\\
\endalign
$$

We focus first on simplifying the term in parenthesis  above.
Since $\lam$ is supported on square-free integers, we may write
$r=\alpha a$, $s=\alpha b$ where $\alpha$, $a$, and $b$ are square-free
with $(a,b)=1$.  Thus $l=\alpha^2 a b$, $l_1=ab$, and $l_2 =\alpha$.
With this notation the sum over $l$ in (6.3) becomes
$$
\align
&\sum_{\alpha} \sum\Sb a, b\\ (a,b)=1\endSb
\frac{r_{\delta}(ab)}{(ab)^{\frac 12 +w}}
\frac{\lam(\alpha a)\lam(\alpha b)}{\alpha^{1+\delta_1+\delta_2}
a^{\frac 12+\delta_1} b^{\frac 12+\delta_2} }
\eta_{\delta}(1+2w;\alpha^2 ab)\\
=& \sum_{\alpha} \frac{1}{\alpha^{1+2\tau}}
\sum\Sb a,b\\ (a,b)=1\endSb \frac{r_{\delta}(a) r_{\delta}(b)}
{a^{1+\delta_1+w}b^{1+\delta_2+w}} \lam(\alpha a) \lam(\alpha b)
\eta_{\delta}(1+2w;\alpha^2 ab).
\tag{6.4}\\
\endalign
$$

Define, for odd primes $p$,
$$
h_{w}(p) = \Big( 1+\frac 1p +\frac 1{p^{1+2w}} -
\frac{p^{-2\delta}+p^{2\delta}}{p^{2+2w}} +\frac{1}{p^{3+4w}}\Big)
$$
and extend this multiplicatively to a function on
odd, square-free integers.  From the definition of $\eta$ we
see that
$$
\eta_{\delta}(1+2w;\alpha^2 ab) =
\frac{\eta_{\delta}(1+2w;1)}{h_{w}(\alpha) h_{w}(a) h_{w}(b)}
\prod_{p|\alpha} \Big( 1+\frac{1}{p^{1+2w}}\Big).
$$
Hence our expression (6.4) may be recast as
$$
\eta_{\delta}(1+2w;1) \sum_{\alpha}
\frac{1}{\alpha^{1+2\tau} h_{w}(\alpha)} \prod_{p|\alpha}\Big(1
+\frac{1}{p^{1+2w}}\Big) \sum\Sb a,b\\ (a,b)=1\endSb
\frac{r_{\delta}(a) \lam(\alpha a)}{a^{1+\delta_1+w}h_{w}(a)}
\frac{r_{\delta}(b)\lam(\alpha b)}{b^{1+\delta_2+w}h_{w}(b)}.
$$
Using $\sum_{\beta|(a,b)} \mu(\beta)= 1$ if $(a,b)=1$ and $0$ otherwise,
the above becomes
$$
\align
\eta_{\delta}(1+2w;1) \sum_{\alpha} \frac{\prod_{p|\alpha}
(1+1/p^{1+2w})}{\alpha^{1+2\tau} h_{w}(\alpha)}
\sum_{\beta} &\frac{r_{\delta}(\beta)^2 \mu(\beta) }{\beta^{2+2\tau
+2w} h_{w}(\beta)^2 } \\
&\times \sum_{a,b}
\frac{r_{\delta}(a) \lam(a\alpha \beta) }{a^{1+\delta_1+w} h_{w}(a)}
\frac{r_{\delta}(b)\lam(b\alpha\beta)}{b^{1+\delta_2+w}h_{w}(b)}.
\tag{6.5}\\
\endalign
$$

Define
for odd primes $p$
$$
H_{w}(p) =  1+ \frac{1}{p^{1+2w}} -\frac{r_{\delta}(p)^2}{
p^{1+2w} h_{w}(p)},
$$
and extend this multiplicatively to all odd, square-free integers.
Then grouping terms according to $\gamma =\alpha \beta$, we
see that (6.6) equals
$$
\eta_{\delta}(1+2w;1) \sum_{\gamma} \frac{H_{w}(\gamma)}
{\gamma^{1+2\tau} h_{w}(\gamma)} \Big( \sum_{a} \frac{r_{\delta}(a)
\lam(a\gamma)}{a^{1+\delta_1+w}h_{w}(a)}\Big)
\Big( \sum_b \frac{r_{\delta}(b) \lam(b\gamma)}{b^{1+\delta_2+w}
h_{w}(b)}\Big). \tag{6.6}
$$

\proclaim{Lemma 6.1}  Let $R$ be a polynomial with $R(0)=R^{\prime}(0)=0$.
Let $g$ be a multiplicative function with $g(p)=1+O(p^{-\nu})$
for some fixed $\nu >0$.  Let $y$ be a large real
number, and suppose that $u$ and $v$ are
bounded complex numbers such that $\text{Re}(u+v)$ and $\text{Re}
(u-v)$ are $\ge - D/\log y$ where $D$ is an absolute positive
constant.  When Re $s>1 + D/\log y$ we have
$$
\sum\Sb n=1\\ n\ \text{odd} \endSb^{\infty}
\frac{r_{v}(n) \mu(nc)}{n^{s+u}} g(n)
= \frac{\mu(c) G(s,c;u,v)}{\zeta(s+u+v)\zeta(s+u-v)},
$$
where $G(s,c;u,v) =\prod_{p} G_p (s,c;u,v)$
with
$$
G_p(s,c;u,v):=
\cases
(1-\frac 1{p^{s+u+v}})^{-1} (1-\frac{1}{p^{s+u-v}})^{-1}
&\text{if } p|2c\\
(1-\frac 1{p^{s+u+v}})^{-1} (1-\frac{1}{p^{s+u-v}})^{-1} (1-\frac{g(p)r_v(p)}
{p^{s+u}})
&\text{otherwise},\\
\endcases
$$
so that $G(s,c;u,v)$ is holomorphic in Re $s> \max(\frac 12,1-\nu)
+ D/\log y$.
For any odd integer $c \le y$ we have
$$
\align
\sum\Sb n\le y/c\\ n \ \text{odd} \endSb
\frac{r_{v}(n)\mu(nc)}{n^{1+u}}&g(n) R\L{\log (y/(cn))}{\log y}
=  O\Big(\frac{E(c)}{\log^2 y} \L{y}{c}^{-\text{Re }u +|\text{Re }v|}
\exp(-A_0\sqrt{\log (y/c)})\Big)\\
&+ \mathop{\text{Res}}_{s=0}
\frac{\mu(c)G(s+1,c;u,v)}{s\zeta(1+s+u+v)\zeta(1+s+u-v)
} \sum_{k=0}^{\infty} \frac{1}{(s\log y)^k} R^{(k)}
\Big(\frac{\log (y/c)}{\log y}\Big)
\\
\endalign
$$
for some absolute constant $A_0>0$, and where $E(c)=\prod_{p|c}
(1+1/\sqrt{p})$.
\endproclaim
\demo{Proof} Our assertion about
the generating function $\sum r_v(n) \mu(nc) g(n)/n^{s+u}$
follows readily upon comparing Euler products.
In proving the other statements we may plainly suppose that $c \le y/2$.
Using the Taylor expansion $R(x)=\sum_{j=0}^{\infty}
\frac{R^{(j)}(0)}{j!} x^j = \sum_{j=2}^{\infty} \frac{R^{(j)}(0)}{j!} x^j$,
we see that our sum is
$$
\align
&\sum_{j=2}^{\infty} \frac{R^{(j)}(0)}{(\log y)^{j}} \frac{1}{j!} \sum\Sb
n\le y/c\\ n\ \text{odd} \endSb \frac{r_{\beta}(n)\mu(nc)}{n^{1+\alpha}} g(n)
\log^{j} \L{y}{cn}
\\
=& \sum_{j=2}^{\infty} \frac{R^{(j)}(0)}{(\log y)^j}
\frac{1}{2\pi i} \int_{(\frac{D+1}{\log (y/c)})}
\frac{\mu(c) G(s+1,c;u,v)}{\zeta(1+s+u+v)\zeta(1+s+u-v)} \Big(\frac yc\Big)^s
\frac{ds}{s^{j+1}}.
\\
\endalign
$$
The integral above is evaluated by a standard procedure:
First one truncates the above integral to the line segment
$\frac{D+1}{\log (y/c)} -iT$ to $\frac{D+1}{\log (y/c)} +iT$ where
$T:=\exp(\sqrt{\log (y/c)})$.  The  error involved in doing so
is $\ll E(c)(\log y/c)^2/T^2$.  Next we shift the integral on this line
segment to the left onto the line segment $-\text{Re }u +|\text{Re } v|
-A_1/\log T$ where
$A_1$ is a positive constant such that $\zeta(1+s+u+v)\zeta(1+s+u-v)$
has no zeros in the region traversed.  We encounter a (multiple)
pole at $s=0$ whose residue we shall calculate presently.  The integrals on
the three other sides are bounded using standard estimates for $1/\zeta(s)$
in the zero-free region, and contribute an amount $\ll E(c)(\log (Ty/c))^2
(T^{-2} + (y/c)^{-\text{Re } u +|\text{Re }v| - A_1/\log T})$.
We conclude that for an
appropriate positive constant $A_0$ the above is
$$
\align
&=\mathop{\text Res}_{s=0}
\frac{\mu(c) G(s+1,c;u,v)}{s\zeta(1+s+u+v)\zeta(1+s+u-v)}
\sum_{j=2}^{\infty}\frac{ R^{(j)}(0) (y/c)^s}{s^{j} \log^j y}
\\
&\hskip 1 in
+ O\Big( \frac{E(c)}{\log^2 y} \L{y}{c}^{-\text{Re }u +|\text{Re }v|}
\exp(-A_0\sqrt{\log (y/c)})\Big).
\\
\endalign
$$
For the purpose of the residue calculation we may replace
$\sum_{j=2}^{\infty}\frac{ R^{(j)}(0) (y/c)^s}{s^{j} \log^j y}$ with
$$
\sum_{j=2}^{\infty}\frac{ R^{(j)}(0)}{s^{j} \log^j y}
\Big(\sum_{l\le j} \frac{s^l}{l!} (\log (y/c))^l\Big)
= \sum_{k=0}^{\infty}
\frac{s^{-k}}{\log^k y}  \sum_{l=0}^{\infty}
\frac{R^{(k+l)}(0)}{l!} \Big(\frac{\log (y/c)}{\log y}\Big)^l,
$$
upon grouping terms according to $k=j-l$, and bearing in mind that
$R(0)=R'(0)=0$.  This clearly equals
$$
\sum_{k=0}^{\infty}
\frac{s^{-k}}{\log^k y} R^{(k)}\Big(\frac{\log (y/c)}{\log y}\Big),
$$
completing our proof of the Lemma.

\enddemo

%Let $P$ be a real-valued polynomial such that $P(0)=P'(0)=0$, and $P(b)=1$,
% $P'(b)=0$ where $b\in (\epsilon,1-\epsilon)$
%is a parameter to be fixed later.
%We shall take $Q(u)=P(u)$ if $u\le b$, and $Q(u)=P(u)+(1-P(u))=1$
%if $b\le u \le 1$.

We now return to the evaluation of
the expression (6.6).
We first deal with the contribution arising from the
terms $\gamma \le M^{1-b}$.  We shall apply Lemma 6.1 twice. In both
cases we take $u=\delta_1+w$, $v=\delta$, and $g(n)=1/h_w(n)$, and
we shall denote the corresponding $G(s,\gamma;u,v)$ by $G_w(s,\gamma;u,v)$.
In the first application we take $y=M$, and $R(x)=P(x)$; and in
the second application we take $y=M^{1-b}$ and $R(x)= (1-P(b+x(1-b)))$.
Adding these two applications we deduce that
$$
\align
\sum\Sb a \le M/\gamma\\ a \ \text{odd} \endSb
\frac{r_{\delta}(a) \mu(a\gamma)}{a^{1+\delta_1+w }
h_w(a)} &Q\Big(\frac{\log(M/a\gamma)}{\log M}\Big)
= \frac{\mu(\gamma)G_w(1,\gamma;\delta_1+w,\delta)}{\zeta(1+\delta_1+w+\delta)
\zeta(1+\delta_1+w-\delta)}\\
&  + O\Big( \frac{E(\gamma)}{\log^2 M} \L{M^{1-b}}
{\gamma}^{-\tau-\text{Re }w} \exp(-A_0 \sqrt{\log M^{1-b}/\gamma})\Big).
\\
\endalign
$$
Observe that, because of our choice of the contour ${\Cal C}$,
the main term above is $\ll |\delta_1|^2$.
An analogous expression holds for the sum over $b$ in (6.6), with
the only change being that $\delta_1$ above gets replaced by $\delta_2$.
We deduce that the contribution of the $\gamma \le M^{1-b}$ terms to
(6.6) equals
$$
\align
\eta_\delta(1+2w;1) &\sum\Sb \gamma \le M^{1-b} \\ \gamma \text{ odd}
\endSb
\frac{H_w(\gamma)}{\gamma^{1+2\tau}h_w(\gamma)}\Big( \frac{\mu^2(\gamma)
G_w(1,\gamma;\delta_1+w,\delta)G_w(1,\gamma;\delta_2+w,\delta)}
{\prod_{\mu=\pm} \zeta(1+\delta_1+w+\mu\delta)\zeta(1+\delta_2+w+\mu\delta)}
\\
& + O\Big( E(\gamma)\frac{|\delta_1|^2 }{\log^2 M}
\L{M^{1-b}}{\gamma}^{-\tau-\text{Re }w} \exp(-A_0 \sqrt{\log (M^{1-b}/\gamma)})
\Big)\Big).\\
\endalign
$$
This is readily seen to be
$$
\align
\eta_\delta(1+2w;1)\sum\Sb \gamma \le M^{1-b} \\ \gamma \text{ odd}
\endSb&
\frac{H_w(\gamma)}{\gamma^{1+2\tau}h_w(\gamma)}\frac{\mu^2(\gamma)
G_w(1,\gamma;\delta_1+w,\delta)G_w(1,\gamma;\delta_2+w,\delta)}
{\prod_{\mu=\pm} \zeta(1+\delta_1+w+\mu\delta)\zeta(1+\delta_2+w+\mu\delta)}
\\
&+ O\Big( \frac{|\delta_1|^2 }{\log^2 M}
M^{(1-b)(-\tau-\text{Re }w)} \Big).
\\
\endalign
$$
We use this expression in (6.3) to evaluate the contribution of the
$\gamma \le M^{1-b}$ terms to the integral there.  From our
choice of ${\Cal C}$, and since $M=X^{\frac 12 -\epsilon}$, we
see that the error term arising from the above is
$$
O\Big( \log^2 X |\delta_1|^3 M^{-2\tau(1-b)}\Big). \tag{6.7}
$$
The main term arising there is
$$
\align
\frac{1}{2\pi i} \int_{{\Cal C}}
&\frac{2{\check\Phi}(w-\tau)}{3\zeta(2)\Gamma_\delta(\tau)}
\L{8X}{\pi}^{w-\tau}
Z(1+2w;\delta) \Gamma_\delta(w) \frac{2w}{w^2-\tau^2}
\\
&\eta_\delta(1+2w;1)\sum\Sb \gamma \le M^{1-b} \\ \gamma \text{ odd}
\endSb \frac{H_w(\gamma)}{\gamma^{1+2\tau}h_w(\gamma)}\frac{\mu^2(\gamma)
G_w(1,\gamma;\delta_1+w,\delta)G_w(1,\gamma;\delta_2+w,\delta)}
{\prod_{\mu=\pm} \zeta(1+\delta_1+w+\mu\delta)\zeta(1+\delta_2+w+\mu\delta)}
dw.\\
\endalign
$$
{\sl A priori} the integrand has two poles (at $\pm \tau$) inside ${\Cal C}$,
but since $\prod_{\mu =\pm}
\zeta(1+\delta_1+w+\mu\delta)^{-1}\zeta(1+\delta_2+w
+\mu\tau)^{-1}$ vanishes (indeed to order $2$)
at $w=-\tau$, in fact we have only
the one simple pole at $w=\tau$.  Thus by Cauchy's theorem
the main term above equals
$$
\frac{2{\check \Phi}(0)}{3\zeta(2)} \frac{\eta_\delta(1+2\tau;1)}
{\zeta(1+2\tau)}
\sum\Sb \gamma \le M^{1-b} \\ \gamma \text{ odd}
\endSb \frac{\mu^2(\gamma)H_\tau(\gamma)}{\gamma^{1+2\tau}h_{\tau}(\gamma)}
G_\tau(1,\gamma;\delta_1+\tau,\delta)G_\tau(1,\gamma;\delta_2+\tau,\delta).
\tag{6.8}
$$

\proclaim{Lemma 6.2}  With $w$ on the contour ${\Cal C}$, and
other notations as above we
have for $x\ge 2$
$$
\align
\eta_\delta(1+2w;1) \sum\Sb \gamma\le x\\ \gamma \ \text{odd} \endSb
&\frac{\mu^2 (\gamma) H_w(\gamma)}{\gamma^{1+2\tau} h_w(\gamma)}
G_w(1,\gamma;\delta_1+w,\delta)G_w(1,\gamma;\delta_2+w,\delta)
\\
&= \zeta(1+2\tau)(1-x^{-2\tau})(1+O(|w-\tau|)) + O(x^{-2\tau}).
\\
\endalign
$$
Further if $1\le y\le x$ then for any smooth function
$R$ on $[0,1]$
$$
\align
\eta_\delta(1+2w;1) \sum\Sb y \le \gamma\le x\\ \gamma \ \text{odd} \endSb
\frac{\mu^2 (\gamma) H_w(\gamma)}{\gamma^{1+2\tau} h_w(\gamma)}
&G_w(1,\gamma;\delta_1+w,\delta)G_w(1,\gamma;\delta_2+w,\delta) R\Big(
\frac{\log \gamma}{\log x}\Big)
\\
&= (1+O(|\delta_1|)) \int_y^x R\Big(\frac{\log t}{\log x}\Big) \frac{dt}
{t^{1+2\tau}}.
\\
\endalign
$$
\endproclaim

\demo{Proof} Upon recalling the definition of $G_w(s,\gamma;u,v)$ from
Lemma 6.1 we see that our desired expression equals
$$
\eta_\delta(1+2w;1)G_{w}(1,1;\delta_1+w,\delta) G_w(1,1;\delta_2+w,\delta)
 \sum\Sb \gamma \le x \\ \gamma \text{ odd}
\endSb \frac{f_w(\gamma)}{\gamma^{1+2\tau}},
$$
say, where $f_w(\gamma)$ is the multiplicative function given by
$$
f_w(\gamma)= \mu^2(\gamma)\frac{H_w(\gamma)}{h_w(\gamma)}
\prod_{p|\gamma} \Big( 1 -\frac{r_\delta(p)}{p^{1+\delta_1+w} h_w(p)}
\Big)^{-1} \Big( 1 -\frac{r_\delta(p)}{p^{1+\delta_2+w} h_w(p)}
\Big)^{-1}.
$$
Plainly $f_w(p)=1+O(1/\sqrt{p})$, say, and so the calculation of
the sum over $\gamma$ becomes a standard exercise.  Writing
the generating function $\sum_{\gamma \ \text{odd}} f_w(\gamma)/\gamma^s
=\zeta(s) F_w(s)$, (note that $F$ is holomorphic in Re $s >\frac 12$), using
Perron's formula, and shifting contours appropriately we deduce easily that
$$
\align
 \sum\Sb \gamma \le x \\ \gamma \text{ odd}
\endSb \frac{f_w(\gamma)}{\gamma^{1+2\tau}} &=
\zeta(1+2\tau) F_w(1+2\tau) - F_w(1)\frac{x^{-2\tau}}{2\tau}
+O( x^{-(\frac 13+2\tau)}) \\
&= \zeta(1+2\tau)(1-x^{-2\tau})F_w(1+2\tau)
+ O(x^{-2\tau}),
\\
\endalign
$$
upon using $\zeta(1+2\tau) = 1/(2\tau) +O(1)$, and that $F_w(1+2\tau) =
F_w(1)+ O(\tau)$.
Denote $\eta_\delta(1+2w;1) G_w(1,1;\delta_1+w,\delta)
G_w(1,1;\delta_2+w,\delta)F_w(1+2\tau)$ by $H(w)$ say.  Then with a little
calculation we may check that $H(w) = H(\tau) + O(|w-\tau|)$, and that
$H(\tau)=1$.  This proves the first assertion of the Lemma.  Our
second claim follows upon using partial summation, and arguing along
similar lines.

\enddemo

Using Lemma 6.2 in (6.8) above, and combining with the
error term estimate (6.7), we conclude that the part of the ${\Cal M}_1(l)$
contribution arising from the $\gamma \le M^{1-b}$ terms equals
$$
\frac{2{\check \Phi}(0)}{3\zeta(2)} ( 1 - M^{-2\tau(1-b)})
+O(\log^2 X |\delta_1|^3 M^{-2\tau(1-b)}).
\tag{6.9}
$$

We now turn to the corresponding contribution from the
terms $\gamma > M^{1-b}$ in (6.6).
Applying Lemma 6.1 with $u$, $v$, $g(n)$, and $G_w(s,\gamma;u,v)$ as above,
and with $R(x)=P(x)$, and $y=M$.   We get that for any odd $M^{1-b} \le
\gamma <M$,  (since $Q(\frac{\log (M/a\gamma)}{\log M})=
P(\frac{\log (M/a\gamma)}{\log M})$
for $\gamma$ in this range)
$$
\align
\sum\Sb a \le M/\gamma\\ a \ \text{odd} \endSb
 &\frac{r_{\delta}(a) \mu(a\gamma)}{a^{1+\delta_1+w }
h_w(a)} Q\Big(\frac{\log(M/a\gamma)}{\log M}\Big)
= O\Big( \frac{E(\gamma)}{\log^2 M }
\L{M}{\gamma}^{-\tau -\text{Re } w} \exp(-A_0\sqrt{\log (M/\gamma)})\Big)
\\
&+\mathop{\text{Res}}_{s=0} \frac{\mu(\gamma) G_w(1+s,\gamma
;\delta_1+w,\delta)}{
s\zeta(1+s+\delta_1+w+\delta)\zeta(1+s+\delta_1+w-\delta)}
\sum_{k=0}^{\infty} \frac{1}{(s\log M)^k} Q^{(k)} \Big(\frac{\log (M/\gamma)}
{\log M}\Big).
\!
\\
\endalign
$$
Write the Taylor expansion of $G_w(1+s,\gamma;u,v)/(\zeta(1+s+\delta_1+w
+\delta)\zeta(1+s+\delta_1+w-\delta))$ as
$a_0+a_1 s+ a_2 s^2 +\ldots$.  Then we see that $a_0= (\delta_1+w+\delta)
(\delta_1+w-\delta)G_w(1,\gamma;\delta_1+w,\delta)
 + O((|\delta_1|+|w|)^3)$, $a_1=2 (\delta_1+w) G_w(1,\gamma;\delta_1+w,\delta)
+ O((|\delta_1|+|w|)^2)$, $a_3= G_w(1,\gamma;\delta_1+w,\delta)
+O(|\delta_1|+|w|)$, and that $a_n \ll_n 1$
for $n \ge 4$.  From this it follows that the residue term above equals
$$
\align
\mu(\gamma)G_w(1,\gamma;\delta_1+w,\delta)
\Big(  (\delta_1+w+\delta)(\delta_1+w-\delta)
&Q\Big(\frac{\log (M/\gamma)}{\log M}\Big) + 2 \frac{\delta_1+w}{\log M}
Q^{\prime}\Big(\frac{\log (M/\gamma)}{\log M}\Big)
\\
&+\frac{1}{\log^2 M}Q^{\prime \prime} \Big(\frac{\log (M/\gamma)}{\log M}\Big)
\Big) + O(|\delta_1|^3).
\\
\endalign
$$
An analogous espression holds for the sum over $b$ in (6.6), replacing
$\delta_1$ above with $\delta_2$.  We use these expressions to
evaluate the contribution to (6.6) from the terms $\gamma > M^{1-b}$.
Firstly, we see that the remainder terms that accrue are (bearing in mind that
$w$ is on the contour ${\Cal C}$)
$$
\align
&\ll \sum_{M^{1-b} \le \gamma \le M}
\frac{1}{\gamma^{1+2\tau}} \Big( \frac{E(\gamma)}{\log^2 M} |\delta_1|^2
\Big(\frac{M}{\gamma}\Big)^{-\tau-\text{Re }w}
\exp(-A_0 \sqrt{\log (M/\gamma)}
) + |\delta_1|^5 \Big)
\\
&\ll \frac{|\delta_1|^2}{\log^2 M} M^{-\tau -\text{Re }w } + M^{-2\tau(1-b)}
|\delta_1|^5 \log M.
\\
\endalign
$$
Secondly, we get that the main term here is (denoting
for brevity $Q^{(j)}(\log (M/\gamma)/\log M)$ by $Q^{(j)}_\gamma$)
$$
\align
\eta_\delta(1+2w;1) \sum\Sb M^{1-b} < \gamma \le M \\
\gamma \ \text{odd} \endSb &\frac{\mu^2(\gamma) H_w(\gamma)}{\gamma^{1+2\tau}
h_w(\gamma)} G_w(1,\gamma;\delta_1+w,\delta)G_w(1,\gamma;\delta_2+w,\gamma)
\\
&\times \Big( (\delta_1+w+\delta)(\delta_1+w-\delta)
Q_\gamma + 2 \frac{\delta_1+w}{\log M}
Q^{\prime}_\gamma
+\frac{1}{\log^2 M}Q^{\prime \prime}_\gamma\Big)\\
&\times
\Big((\delta_2+w+\delta)(\delta_2+w-\delta)
Q_\gamma + 2 \frac{\delta_2+w}{\log M}
Q^{\prime}_\gamma
+\frac{1}{\log^2 M}Q^{\prime \prime}_\gamma\Big).\\
\endalign
$$
Applying Lemma 6.2 (and a suitable change of variables)
we conclude that this equals
$$
\align
\log M \int_0^b M^{-2\tau(1-x)}
&\Big( (\delta_1+w+\delta)(\delta_1+w-\delta)
Q(x) + 2 \frac{\delta_1+w}{\log M}
Q^{\prime}(x)
+\frac{Q^{\prime \prime}(x)}{\log^2 M}\Big)
\\
&\times \Big((\delta_2+w+\delta)(\delta_2+w-\delta)
Q(x) + 2 \frac{\delta_2+w}{\log M}
Q^{\prime}(x)
+\frac{Q^{\prime\prime}(x)}{\log^2 M}\Big) dx \\
&\hskip 1 in + O(M^{-2\tau(1-b)} |\delta_1|^5
\log M) . \tag{6.10}
\\
\endalign
$$

We use these expressions in (6.3) to evaluate the
contribution of the $\gamma >M^{1-b}$ terms to the integral there.
From our choices of $M$ and ${\Cal C}$ we get that the error term arising
from the above is
$$
O\Big( |\delta_1|^3  M^{-2\tau} \log^2 X+ |\delta_1|^6
 M^{-2\tau(1-b) }\log^5 X\Big). \tag{6.11}
$$
Call the main term in (6.10) $N(w)$.  Inserting this into the integral in
(6.3), we seek to evaluate
$$
\frac{1}{2\pi i} \int_{\Cal C} \frac{2{\check \Phi}(w-\tau)}{3\zeta(2)
\Gamma_\delta(\tau)}
\L{8X}{\pi}^{w-\tau} Z(1+2w;\delta) \Gamma_\delta(w) N(w) \frac{2w}
{w^2-\tau^2}dw.
$$
Now ${\check \Phi}(w-\tau)\Gamma_\delta(w)/\Gamma_\delta(\tau) =
{\check \Phi}(0) +O(|\delta_1|)$, and $2w Z(1+2w;\delta)=
\frac{1}{4(w^2-\delta^2)} + O( |\delta_1| \log^2 X)$, and
$N(w) \ll M^{-2\tau(1-b)}|\delta_1|^4 \log M$, whence we deduce that
the above integral is
$$
\align
\frac{2{\check \Phi}(0)}{3\zeta(2)} \frac{1}{2\pi i} \int_{\Cal C}
&\L{8X}{\pi}^{w-\tau} N(w) \frac{1}{4(w^2-\delta^2)} \frac{1}{w^2-\tau^2}
dw +O(|\delta_1|^6 M^{-2\tau(1-b)}\log^5 X)\\
=&\frac{2{\check \Phi}(0)}{3\zeta(2)} \frac{1}{8\delta_1 \delta_2 \tau}
\Big( N(\tau) - \L{8X}{\pi}^{-2\tau} N(-\tau)\Big) +
O(|\delta_1|^6 M^{-2\tau(1-b)}\log^5 X). \tag{6.12}\\
\endalign
$$
Using integration by parts together with $Q(0)=Q^{\prime}(0)=0$,
and $Q(b)=1$, $Q^{\prime}(b)=0$ we may simplify the expression
for $N(\tau)$ to
$$
8\delta_1\delta_2 \tau M^{-2\tau(1-b)}
+ \frac{4\delta_1 \delta_2}{\log M}
\int_0^b M^{-2\tau(1-x) } \Big|
Q^{\prime}(x) +\frac{Q^{\prime\prime}(x)}{2\delta_1 \log M}
\Big|^2
%\Big( Q^{\prime}(x)+\frac{Q^{\prime\prime}(x)}{2\delta_2 \log M}
%\Big)
dx.
$$
Similarly we find that
$$
N(-\tau) = \frac{4\delta_1 \delta_2}{\log M}
\int_0^b M^{-2\tau(1-x) } \Big|
Q^{\prime}(x) +\frac{Q^{\prime\prime}(x)}{2\delta_1 \log M}
\Big|^2
%Big( Q^{\prime}(x)+\frac{Q^{\prime\prime}(x)}{2\delta_2 \log M}
%\Big)
dx.
$$
Using these identities we conclude that our expression in (6.12)
equals
$$
\align
\frac{2{\check \Phi}(0)}{3\zeta(2)} \Big( M^{-2\tau(1-b)}
+ \frac{1-(8X/\pi)^{-2\tau} }{2\tau\log M}
&\int_0^b M^{-2\tau(1-x) } \Big|
Q^{\prime}(x) +\frac{Q^{\prime\prime}(x)}{2\delta_1 \log M}
\Big|^2
% \Big( Q^{\prime}(x)+\frac{Q^{\prime\prime}(x)}{2\delta_2 \log M}
%\Big)
dx \Big)
\\
&+O(M^{-2\tau(1-b)}
|\delta_1|^6 \log^5 X).\\
\endalign
$$
Combining this with (6.11) we conclude that the part of
the ${\Cal M}_1(l)$ contribution arising from the $M^{1-b}
\le \gamma \le M$ terms equals
$$
\align
\frac{2{\check \Phi}(0)}{3\zeta(2)} \Big( M^{-2\tau(1-b)}
&+ \frac{1-(8X/\pi)^{-2\tau} }{2\tau\log M}
\int_0^b M^{-2\tau(1-x) } \Big|
Q^{\prime}(x) +\frac{Q^{\prime\prime}(x)}{2\delta_1 \log M}
\Big|^2
% \Big( Q^{\prime}(x)+\frac{Q^{\prime\prime}(x)}{2\delta_2 \log M}
%\Big)
dx \Big)
\\
&+O(|\delta_1|^3 M^{-2\tau} \log^2 X +
M^{-2\tau(1-b)} |\delta_1|^6 \log^5 X)
\endalign
$$

Taking this together with (6.9) we have determined the ${\Cal M_1}(l)$
contribution to be
$$
\align
\frac{2{\check \Phi}(0)}{3\zeta(2)} \Big( 1+
\frac{1-(8X/\pi)^{-2\tau} }{2\tau\log M}  &\int_0^b M^{-2\tau(1-x) } \Big|
Q^{\prime}(x) +\frac{Q^{\prime\prime}(x)}{2\delta_1 \log M}
\Big|^2
%\Big( Q^{\prime}(x)+\frac{Q^{\prime\prime}(x)}{2\delta_2 \log M}
%\Big)
dx \Big)\\
&+O(M^{-2\tau(1-b)} |\delta_1|^6 \log^5 X).
\tag{6.13} \\
\endalign
$$

The calculation of the ${\Cal M}_2(l)$ contribution to (6.2a) is
entirely similar.  We obtain that this contribution equals
$$
\align
-\frac{2{\check \Phi}(0)}{3\zeta(2)} \L{8X}{\pi}^{-\tau}
\frac{(8X/\pi)^{\delta}-(8X/\pi)^{-\delta}}{2\delta \log M}
& \int_{0}^{b} M^{-2\tau(1-x) } \Big|
Q^{\prime}(x) +\frac{Q^{\prime\prime}(x)}{2\delta_1 \log M}
\Big|^2
% \Big( Q^{\prime}(x)+\frac{Q^{\prime\prime}(x)}{2\delta_2 \log M}
%\Big)
dx\\
& O(X^{-\tau} M^{-2\tau(1-b)} |\delta_1|^6 \log^5 X).
\\
\endalign
$$
Inputing this and (6.13) into (6.2a), we
obtain Proposition 2.4.
%we conclude at last that
%$$
%\align
%{\Cal S}(|L(\tfrac 12+\delta_1,\chi_{-8d})&M(\tfrac 12+\delta_1,
%d)|^2; \Phi) = O(M^{-2\tau(1-b)} |\delta_1|^6 \log^5 X + X^{-\tau-\epsilon}
%\Phi_{(2)}\Phi_{(3)}^{\epsilon})\\
%&+\frac{2{\check \Phi}(0)}{3 \zeta(2)}
%\Big(1 + \Big( \frac{1-(8X/\pi)^{-\tau}}{2\tau \log M}
%- \L{8X}{\pi}^{-\tau} \frac{(8X/\pi)^{\delta}-(8X/\pi)^{-\delta}}
%{2\delta\log M}\Big)\\
%&\hskip 1 in \times  \int_{0}^{b} M^{-2\tau(1-x) } \Big|
%Q^{\prime}(x) +\frac{Q^{\prime\prime}(x)}{2\delta_1 \log M}
%\Big|^2 dx\Big).\tag{5.20}\\
%\endalign
%$$

\enddocument

\end